\documentclass[preprint,12pt]{elsarticle}

\usepackage{geometry}
\usepackage[pagebackref=true]{hyperref}

\usepackage{amssymb,amsthm,amsmath}
\usepackage{latexsym}
\usepackage{mathrsfs}
\usepackage{stmaryrd}
\usepackage{graphicx}
\usepackage{color}
\usepackage{enumerate,dsfont}
\usepackage{ushort,phonetic,cancel}

\newtheorem{theorem}{Theorem}[section]
\newtheorem{lemma}[theorem]{Lemma}
\newtheorem{prop}[theorem]{Proposition}

\theoremstyle{definition}
\newtheorem{definition}[theorem]{Definition}

\theoremstyle{remark}
\newtheorem{remark}[theorem]{Remark}
\newtheorem{corr}[theorem]{Corollary}

\newcommand{\C}{ \mathbb{C}}
\newcommand{\R}{\mathbb{R}}
\newcommand{\Z}{ \mathbb{Z}}
\newcommand{\N}{ \mathbb{N}}

\newcommand{\jj}{{\rm i}}
\newcommand{\ee}{{\rm e}}
\newcommand{\dd}{{\rm d}}

\newcommand{\Lop}{{\rm L}}
\newcommand{\Dop}{{\rm D}}
\newcommand{\Iop}{{\rm I}}

\newcommand{\CF}{\widehat{\mathscr{P}}}

\newcommand{\ab}{{b\hspace{-0.45em}a}}
\def\CC{{ \cancelto{}{C} }}
\def\SS{{ \mathcal{S}}}

\DeclareMathOperator*{\real}{Re}
\DeclareMathOperator*{\imag}{Im}

\newcommand{\indc}{{\mathds{1}}}
\newcommand{\Dord}{{\gamma}}

\graphicspath{{./Figures/}} 

\journal{X}

\begin{document}

\begin{frontmatter}



\title{Complex-Order Scale-Invariant Operators and Self-Similar Processes}


\author[Sharif]{Arash Amini}
\author[EPFL2]{Julien Fageot}
\author[EPFL]{Michael Unser}

\address[Sharif]{Advanced Communication Research Institute (ACRI), EE department, Sharif University of Technology, Tehran, Iran.}

\address[EPFL2]{AudioVisual Communications Laboratory (LCAV), \'Ecole polytechnique f\'ed\'erale de Lausanne (EPFL), Lausanne, Switzerland.}

\address[EPFL]{Biomedical Imaging Group (BIG), \'Ecole polytechnique f\'ed\'erale de Lausanne (EPFL), Lausanne, Switzerland.}

\begin{abstract}
Derivatives and integration operators are well-studied examples of linear operators that commute with scaling up to a fixed multiplicative factor; i.e., they are scale-invariant. Fractional order derivatives (integration operators) also belong to this family. In this paper, we extend the fractional operators  to complex-order operators by constructing them in the Fourier domain. We analyze these operators in details with a special emphasis on the decay properties of the outputs. We further use these operators to introduce a family of complex-valued stable processes that are self-similar with complex-valued Hurst indices. These processes are expressed via the characteristic functionals over the Schwartz space of functions. Besides the self-similarity and stationarity, we  study the regularity (in terms of Sobolev spaces) of the proposed processes.  

\end{abstract}


\end{frontmatter}


\section{Introduction}

    \subsection{Scale-invariance Operators and Self-similar Random Processes}
    
Linear differential operators are of great interest to model physical phenomena \cite{Lanczos2012}. They are  tightly linked with mathematical notions such as splines \cite{Unser2007} and stochastic processes \cite{Blu2007}. The derivative operator $\Dop=\tfrac{\dd}{\dd x}$ is the simplest differential operator with a non-trivial null-space. Beside being linear, this operator is also shift- and scale-invariant. 
In general, any operator $\Lop$ is called scale-invariant if it commutes with scaling up to some proportional factor. Formally, this means that
${\Lop  \{ f (T \cdot ) \} = c_T \Lop \{f \} (T \cdot) }$, 
where $c_{T} > 0  $ is a multiplicative scalar that depends on $T > 0$ and not on $f$. It is known that $c_{T}=T^{\Dord}$ for some fixed $\Dord$ that depends on $\Lop$ \cite{Gelfand_Shilov}, meaning that
\begin{align} \label{eq:scaleinvdefinition}
{\Lop  \{ f (T \cdot ) \} = T^\gamma \Lop \{f \} (T \cdot) } 
\end{align}
for every $f$ and $T$. 
For instance, when $\Lop = \Dop$, we have  that $\Dord = 1$.
One can easily verify that all $\Dop^k$ operators ($k\in\N$) are scale-invariant with $\Dord=k$.
More generally, the Riemann-Liouville fractional derivatives $\Dop^{\Dord}$ with $\gamma\in\R^{+}$ are also scale-invariant. However,
$\Dop^{\gamma}$ is not local when $\gamma$ is not an integer.

A real-valued random process $S = (S(x))_{x \in \R}$ is called \textit{self-similar} 
if, {for each $T > 0$, there exists a constant $d_T > 0$ such that}
${S( T \cdot )  \stackrel{d}{=} d_T S.}$
This means that the random processes $S(T \cdot) = (S(Tx))_{x\in \R}$ and $d_T S = (d_T S(x))_{x\in \R}$ have the same finite-dimensional marginal laws for every $T>0$. It is known that $d_T = T^H$ for some $H \in \R$ called the Hurst index of $S$~\cite{Embrechts2000}. This implies that
\begin{align}\label{eq:SSprocDef}
{S( T \cdot )  \stackrel{d}{=} T^H S.}
\end{align}
for any $T$. 
The fBm processes cover the whole range of $0<H<1$, with $H=\tfrac{1}{2}$ representing the standard Brownian motion. 
Moreover, the fractional Brownian motion (fBm) processes are the only self-similar Gaussian processes with stationary increments \cite{SamorodnitskyBook94}.
While the fBm process with Hurst index $0<H<1$ is not stationary, it can be whitened via the fractional-derivative operator of order $\Dord= H+\tfrac{1}{2}$ \cite{UnserBook2014}. 


\subsection{Related Works}

The specification of fBm by Mandelbrot and van Ness\footnote{The fBm processes were previously studied by Kolmogorov \cite{Kolmogorov1940} and L\'evy \cite{Levy1953}; however, they were made popular by Mandelbrot and van Ness.} in \cite{Mandelbrot1968}  led to a substantial increase in the popularity of fractals and fractional derivatives. 
Mandelbrot and van Ness define fBM as
\begin{align} \label{eq:fBm}
B_H(x) = \frac{1}{\Gamma(H+\frac{1}{2})} \int_{\R} \Big( (x-\tau)_{+}^{H-\frac{1}{2}} - (-\tau)_{+}^{H-\frac{1}{2}} \Big) \, \dd B(\tau), 
\end{align}
where the random process $B(x)$ is the classical Brownian motion and the notation $(x)_{\pm}^{r}$ refers to $|x|^r \indc_{\pm x>0}$ with $\pm$ representing either $+$ or $-$, but consistent in both. Due to their self-similarity and long-range-dependence properties, fBm found applications in various fields such as traffic modeling \cite{Leland1994}, image processing \cite{Pentland1984,Lundahl1986,Chen1989,pesquet2002stochastic}, modeling scatterings from rough surfaces \cite{Franceschetti1999}, and finance \cite{Rogers1997}. The self-similarity property implies that these random processes have the same structure at various scales. As wavelets provide a multiresolution representation of signals, the study of the wavelet representation of fBm processes has been the center of extensive research \cite{Flandrin1992,Meyer1999, Hwang1994, Pipiras2005, Tafti2009, Fageot2015, Storath2017,Fageot2017nterm}.

The extension of fBm  has been considered in a number of aspects. While early extensions generalized the random processes to two \cite{Heneghan1996} and higher dimensions \cite{Reed1995}, the introduction of non-Gaussian distributions were another generalization. Indeed, Gaussian distributions belong to the larger family of $\alpha$-stable distributions ($0<\alpha\leq 2$). As shown in \cite{SamorodnitskyBook94}, the random processes
\begin{align}\label{eq:FracStableMotion}
S_{\alpha,H}(a,b;x) = \int_{\R} \Big( & a\big( (x-\tau)_{+}^{H-\frac{1}{\alpha}} - (-\tau)_{+}^{H-\frac{1}{\alpha}} \big) \nonumber\\
&+ b\big( (x-\tau)_{-}^{H-\frac{1}{\alpha}} - (-\tau)_{-}^{H-\frac{1}{\alpha}} \big)\Big) \, \dd \mu_{\alpha}(\tau),
\end{align}
are also self-similar with $\alpha$-stable marginal distributions, where $a,b$ are arbitrary reals, $H$ is any real in the interval $]0,1[$ other than $\frac{1}{\alpha}$, and $\mu_{\alpha}$ is a suitable $\alpha$-stable motion. 
The case of $H=\frac{1}{\alpha}$ could be defined as the conventional stable motion process, while $H=1$ represents a  degenerate case (lines with random slopes). 
Indeed, $S_{2,H}\big(\frac{1}{\Gamma(H+\frac{1}{2})},0;x\big)$ associated with $\alpha=2$ is equivalent to \eqref{eq:fBm}. 
This wider class of distributions allows for a wider range of applications \cite{Watkins2009, Burnecki2010} and with better modeling capabilities \cite{Mikosch2002}. 
For $\alpha \neq 2$, the stable distributions have heavy tails with infinite variance.
Self-similar stable processes have been extensively studied~\cite{SamorodnitskyBook94,kono1991self,Fageot2019scaling,Fageot2017besov}.


The stationary CARMA processes are usually expressed via stochastic integrals $L_h(x)= \int h(x-\tau) \, \dd L(\tau)$ where $L$ is a L\'evy process and $h$ is the kernel associated with the CARMA differential equation~\cite{Brockwell2000CARMA,Brockwell2009existence}.
By contrast, their fractional variants (FICARMA processes) make  explicit use of the fractional derivative/integration operators  \cite{Brockwell2005,MarquardtThesis}. These are defined by replacing the kernel $h$ with its fractionally integrated version $h_{\Dord}(x) = \big(I_{+}^{\Dord} h\big)(x)$, where
\begin{align}\label{eq:FICARMA}
\big(I_{\pm}^{\Dord} \varphi\big)(x) = \frac{1}{\Gamma(\Dord)} \int_{\R} \varphi(\tau)\, (x-\tau)_{\pm}^{\Dord} \dd \tau.
\end{align}
Thus, FICARMA processes are stationary and can be whitened by applying the suitable linear shift-invariant operator. If $h(\cdot)=\delta(\cdot)$ was admissible mathematically in \eqref{eq:FICARMA} (this is not the case in \cite{Brockwell2005,MarquardtThesis}), then, this fractional extension would yield 
stationary fractional L\'evy processes. 

The same idea is developed in \cite{Huang2007} with a different technique. 
Based on an earlier interpretation of fBm as functionals of the Gaussian white noise in \cite{Huang2006}, 
fractional $\alpha$-stable motions (and sheets) are defined in \cite{Huang2007} for 
 $\alpha\in]1,2[$ and {$\Dord\in ]1,2-\tfrac{1}{\alpha}[$. 
 They are shown to be self-similar with Hurst index $H = \Dord + \frac{1}{\alpha} - 1 \in ]  \frac{1}{\alpha} , 1[$~\cite[Proposition 4.2]{Huang2006}. }
The random processes are defined as generalized random processes; \emph{i.e.}, via their probability law over the space of tempered generalized functions~\cite{Gelfand_Vilenkin,Ito1984foundations}. 
This is also the approach that we shall pursue in this paper.

While the range of the Hurst exponent $H$ is dominantly limited to $]0,1[$ in most research works devoted to fBm processes, there have been few attempts to enlarge this range. The existence of Gaussian processes 
\begin{align}
B_H(x)=\frac{1}{\Gamma(H+\frac{1}{2})} \int_{\R} \Big( (x-\tau)_{+}^{H-\frac{1}{2}}  -\sum_{j=0}^{k-1} \frac{x^j}{j!}  \tfrac{\dd^j}{\dd \tau^{j}} (-\tau)_{+}^{H-\frac{1}{2}}\Big) \dd B(\tau)
\end{align}
with $H\in]k-1,k[$, $k \in \mathbb{N}$, is briefly mentioned in \cite{Bel1998}; these random processes have stationary $k$th order increments. More details about these random processes with $H$ covering all non-integer positive values is provided in \cite{Perrin2001}. Another more systematic point of view is provided in \cite{TaftiThesis}, where
fractional processes are  defined as the solution to fractional stochastic differential equations. 
The solutions are constructed by incorporating the adjoint operator of adequate fractional integration operators in the characteristic functional of the random process. One of the advantages of this approach is that the operators and the canonical {probability} measures can be studied separately and then combined to form the desired processes. 
The properties of the involved fractional operators and their adjoints are investigated in \cite{Sun2012}. Characterizing such random processes as solutions of stochastic differential equations provides a unified framework for the stationary CARMA, non-stationary fractional L\'evy, and some higher-dimensional random fields {\cite{UnserBook2014,Unser2014_I, Unser2014_II,Fageot2014,fageot2017gaussian}.}

\subsection{Contributions and Outline}

The connection between self-similar processes, operators, and splines in \cite{Unser2007} and, particularly, complex-order B-splines in \cite{Forster2006} and complex-order exponential splines in \cite{Massopust2013}, are the inspirations for this work.  
Specifically, our goal is to construct self-similar random processes with \textit{complex-valued} Hurst index $H$ with the help of scale-invariant operators with \textit{complex-valued} order  that can be whitened.
Our contributions can be summarized as follows.


\begin{enumerate}
\item \textit{Systematic characterization of scale-invariant operators with complex order.}
We introduce an extended family of scale-invariant derivatives and integration operators. Their respective construction is made in Fourier domain, and allows for complex-valued $\gamma$ in \eqref{eq:scaleinvdefinition}. 
We provide an in-depth analysis of these fractional operators, with a special emphasis on the decay properties of the outputs. 
We believe these results to be interesting by themselves, but they are also preparatory to the construction of self-similar processes. 

\item \textit{Construction of novel self-similar stable processes with complex Hurst exponent.}
We introduce an extended family of fractional stable processes  that are self-similar as in \eqref{eq:SSprocDef}, but with complex-valued Hurst index $H$. Obviously, a complex-valued $H$ implies that the random process itself  is complex-valued as well. 
Similarly to \cite{UnserBook2014}, we specify the random processes via characteristic functionals over the (complex-valued) Schwartz space $\mathcal{S}$ of rapidly decaying and smooth test functions and invoke the Bochner-Minlos theorem to conclude the existence of a probability measure over $\mathcal{S}^{'}$, the space of (complex-valued) tempered generalized functions.  
%
We moreover  study the invariance (self-similarity and stationarity) and regularity (in terms of Sobolev spaces) of the proposed fractional stable random processes. \\
\end{enumerate}

The paper is organized as follows. Complex order fractional derivative and integral operators are introduced and analyzed in Section~\ref{sec:Theos}. The family of self-similar stable processes is defined and studied in Section~\ref{sec:processes}. Section~\ref{sec:Lemmas} collects useful results for the proofs of the main results, which are provided in Section~\ref{sec:Proofs}. Finally, we conclude in Section~\ref{sec:conclusion}.


\section{Complex Fractional Operators and their Properties}\label{sec:Theos}

Our two classes of complex-order fractional  operators --derivatives vs. intregrators-- are introduced in Section~\ref{sec:constructop}. 
The properties of derivative operators are studied in Section~\ref{sec:propertiesderivatives}.
Since the integration operators, unlike the derivative operators, are not necessarily shift-invariant, we separately study them in Section~\ref{sec:propertiesintegration}.


\subsection{Construction of Complex-Order Fractional Operators} \label{sec:constructop} 

The Fourier and inverse Fourier transforms, $\mathcal{F}\big\{\cdot\big\}$ and $\mathcal{F}^{-1}\big\{\cdot\big\}$, respectively, are understood  as
\begin{align}
\widehat{f}(\omega) &= \mathcal{F}\big\{f\big\}(w) = \int_{\R} f(x) \ee^{-\jj \omega x} \dd x, \\
f(x) &= \mathcal{F}^{-1}\big\{ \widehat{f}\big\}(x) = \tfrac{1}{2\pi}\int_{\R} \widehat{f}(\omega) \ee^{\jj \omega x} \dd \omega,
\end{align}
for any $x, \omega \in \R$,
{provided that $f$ and/or $\widehat{f}$ are integrable.}
The circumflex accent ($\,\widehat{\cdot}\,$) always denotes the Fourier transform of a function with the symbol 
$\omega$ used for the frequency variable. 
It is well-known that the (inverse) Fourier transforms of Schwartz functions are themselves Schwartz functions. Further, all {tempered generalized functions} also have valid (inverse) Fourier transforms in the distributional sense~\cite{Schwartz1966distributions}. The Fourier transforms of $f(Tx)$ and of $\tfrac{\dd}{\dd x} f(x)$ are given by $\tfrac{1}{|T|}\widehat{f}(\tfrac{\omega}{T})$ and $\jj \omega \widehat{f}(\omega)$, respectively. Accordingly, 
all integer-order derivative operators $\Dop^{k}$ can be redefined as
\begin{align}\label{eq:Integer_Deriv_Fourier}
\big(\Dop^{k} \varphi\big)(x) = \mathcal{F}^{-1}\big\{ (\jj \cdot)^k \widehat{\varphi}(\cdot) \big\}(x).
\end{align}
To extend the framework to fractional exponents, we define
\begin{align}\label{eq:hDef}
h_{a,b}^{\Dord}(\omega) = \left\{\begin{array}{cl}
a\, (\omega)_{+}^{\Dord} + b\, (\omega)_{-}^{\Dord}, & \omega \neq 0, \phantom{\Big|}\\
0, & \omega = 0, \phantom{\Big|}
\end{array}
\right.
\end{align}
where $a,b \in \C \backslash\{0\}$ and  $\Dord \in \C$.
It is easy to check that $h_{a,b}^{\Dord}(T\omega) = T^{\Dord} h_{a,b}^{\Dord}(\omega)$ for all $T > 0$ and $\omega\in \R$. In other words, the dilation of $h_{a,b}^{\Dord}$ translates into an amplitude scaling. The generalized functions that satisfy this property are called homogeneous. In fact, except for singularities at $\omega=0$ (Dirac's impulse and its derivatives),  \eqref{eq:hDef} describes the complete family of  homogeneous generalized functions \cite{Gelfand_Shilov}. It follows that any well-behaved\footnote{As before, we are excluding point singularities such as $\delta(\omega)$ and its derivatives.} linear and shift-invariant operator $\Lop$ (acting on $\mathcal{S}$) that is also scale-invariant corresponds to a Fourier multiplier of the form  $h_{a,b}^{\Dord}$. By comparing these operators with \eqref{eq:Integer_Deriv_Fourier}, we observe that such scale-invariant operators can be interpreted as $\Dord$-order derivatives. 

\begin{definition}[\bf Complex-order derivatives] \label{def:Dabgamma}
The $\gamma$-order derivative operator  $\Dop_{a,b}^{\Dord}$ with parameters $\gamma \in \C$ with $\real(\gamma) > -1$ and $a, b \in \C \backslash\{0\}$ is defined as 
$\Dop_{a,b}^{\Dord}$ as
\begin{align}\label{eq:fracD_Def}
\big(\Dop_{a,b}^{\Dord} \, \varphi\big)(x) 
= \mathcal{F}^{-1}\Big\{\widehat{\varphi}(\cdot) h_{a,b}^{\Dord}(\cdot) \Big\}(x)
= \tfrac{1}{2\pi} \int_{\R} \widehat{\varphi}(\omega) h_{a,b}^{\Dord}(\omega) \ee^{\jj \omega x} \dd \omega
\end{align}
for any $x \in \R$ and $\varphi \in \mathcal{S}$. 
\end{definition}

{Definition \ref{def:Dabgamma} is valid because the function $h_{a,b}^\gamma$ specifies a tempered generalized function: it is integrable at the origin for $\real(\gamma) > -1$ and asymptotically grows like a polynomial. It is therefore the frequency response; \emph{i.e.} the Fourier transform of the impulse response of a convolution operator $\Dop_{a,b}^{\Dord} : \mathcal{S} \mapsto \mathcal{S}'$ 
\cite{Schwartz1966distributions}.}

{Classical fractional operators are covered by Definition~\ref{def:Dabgamma} when restricting to $\gamma \in \R$. The fractional Laplacian $(-\Delta)^{\gamma/2}$, whose Fourier multiplier is $\omega \mapsto |\omega|^\gamma$, corresponds to $a=b=1$ in \eqref{eq:hDef}. The fractional derivative $\Dop^\gamma$ is obtained with $a = \mathrm{i}^\gamma$ and $b = (-\mathrm{i})^\gamma$.}

For $\real (\Dord)>0$, the Fourier multiplier $h_{a,b}^{\Dord}(\omega)$ associated with fractional derivatives vanishes at $\omega=0$. Hence, the operator might not be invertible (similar to the integer-order case). 
However, we can introduce integration operators as possible inverses of the fractional derivatives as follows. We denote by $\N$ the set of integers equal or larger than $1$.

{\begin{definition}[\bf Complex-order integrators] \label{def:Igammakab}
The $\gamma$-order integration operator $\Iop_{a,b}^{\Dord; \scalebox{0.6}{$k$}}$ with parameters $\gamma \in \C$ with $\real(\gamma) > 0$, $a, b \in \C \backslash\{0\}$, and $k \in \mathbb{N}$ such that $k\geq \lfloor \real(\Dord) \rfloor$ is defined as 
\begin{align}\label{eq:fracI_Def}
\big(\Iop_{a,b}^{\Dord; \scalebox{0.6}{$k$}}\, \varphi\big)(x) 
=\tfrac{1}{2\pi} \int_{\R} \widehat{\varphi}(\omega) \frac{\ee^{\jj \omega x} - \sum_{j=0}^{k-1}\frac{(\jj x)^j}{j!}\omega^j }{h_{a,b}^{\Dord}(\omega)}
  \dd \omega,
\end{align}
for every $x \in \R$ and $\varphi \in \mathcal{S}$. 
\end{definition}}

{Again, one can verify that Definition \ref{def:Igammakab} is valid since the function $\omega \mapsto \frac{\ee^{\jj \omega x} - \sum_{j=0}^{k-1}\frac{(\jj x)^j}{j!}\omega^j }{h_{a,b}^{\Dord}(\omega)}$ is locally integrable (in particular around the origin, since $k \geq \lfloor \real(\gamma) \rfloor$) and bounded by some polynomial. Hence $\big(\Iop_{a,b}^{\Dord; \scalebox{0.6}{$k$}}\, \varphi\big)$ specifies a tempered generalized function.}

More precisely, we shall show that $\Iop_{a,b}^{\Dord; \scalebox{0.6}{$k$}}$ forms a proper right-inverse for $\Dop_{a,b}^{\Dord}$ if $k=\lfloor \real(\Dord)\rfloor$ or $k= \lceil \real(\Dord) \rceil$ (see Proposition~\ref{prop:Integration_Adjoint}). Moreover, the adjoint operators $\Iop_{a,b}^{(\Dord; \scalebox{0.6}{$k$})*}$ 
for such integration operators are characterized by the relation 
$\langle  \Iop_{a,b}^{(\Dord; \scalebox{0.6}{$k$})*} \varphi , \psi \rangle = \langle \varphi , \Iop_{a,b}^{\Dord; \scalebox{0.6}{$k$}} \psi \rangle$ for any $\varphi, \psi \in \mathcal{S}$ and 
are given by
\begin{align}\label{eq:fracIAdj_Def}
\big(\Iop_{a,b}^{(\Dord; \scalebox{0.6}{$k$})*}\, \varphi\big)(x) 
= \mathcal{F}^{-1}\left\{ \frac{ \widehat{\varphi}(\cdot)-\sum_{j=0}^{k-1}\frac{ \widehat{\varphi}^{(j)}(0) }{j!}(\cdot)^j }{h_{a,b}^{\Dord}(-\cdot)} \right\}(x).
\end{align}
Next, we characterize the key properties of the fractional derivative and integration operators $\Dop_{a,b}^{\Dord}$, $\Iop_{a,b}^{\Dord; \scalebox{0.6}{$k$}}$ and  $\Iop_{a,b}^{(\Dord; \scalebox{0.6}{$k$})*}$.

\subsection{Properties Fractional derivative operators}
\label{sec:propertiesderivatives} 
Although fractional derivatives are well-documented in the literature \cite{samko1993fractional}, we are not aware of any systematic study of 
complex order operators. The first issue is to characterize  the smoothness and decay properties of $\Dop_{a,b}^{\Dord} \,\varphi$ for a given test function $\varphi$.

\begin{theorem}\label{theo:decayRate}
Let $a,b,\Dord$ be arbitrary complex numbers with $\real(\Dord) >-1$, and let $h_{a,b}^{\Dord}$ be the homogeneous generalized function of degree $\Dord$ defined in \eqref{eq:hDef} corresponding to the $\Dord$-order derivative operator $\Dop_{a,b}^{\Dord}$ (defined in \eqref{eq:fracD_Def}). Then, $\big(\Dop_{a,b}^{\Dord} \varphi\big)(x)$ is infinitely differentiable for $\varphi\in\SS$. Further,
\begin{enumerate}[(i)] 
\item \label{decayTheo_claim:gooddecay}
if $\real(\Dord) \not \in {\N}\cup\{0\}$, then, $\Dop_{a,b}^{\Dord}$ is a continuous mapping from $\SS$ into $L_p$ for $p>\frac{1}{\real(\Dord) +1}$, and there exists a constant $c > 0$ such that
\begin{align*}
\forall \varphi \in\SS,\,\, x\in \R:~~~ \Big| \big(\Dop_{a,b}^{\Dord} \varphi\big)(x)\Big| \leq \frac{c}{1 + |x|^{\real(\Dord) +1}}.
\end{align*}

\item \label{decayTheo_claim:schwarz}
if $\Dord \in {\N}$ and $h_{a,b}^{\Dord}(\omega) \equiv d \omega^{\Dord}$ for some $d\in\C$, then, $\Dop_{a,b}^{\Dord}$ is a continuous mapping from $\SS$ into $\SS$.

\item \label{decayTheo_claim:slowdecay}
if $\real(\Dord) \in {\N}$ but $h_{a,b}^{\Dord}(\omega) \not\equiv c \omega^{\Dord}$, then, $\Dop_{a,b}^{\Dord}$ is a continuous mapping from $\SS$ into $L_p$ for $p>\frac{1}{\real(\Dord) +1}$, and there exists a constant $c > 0$ such that
\begin{align*}
\forall \varphi \in\SS,\,\, x\in \R:~~~ \Big| \big(\Dop_{a,b}^{\Dord} \varphi\big)(x)\Big| \leq c\frac{1+\log(1+|x|)}{1 + |x|^{\real(\Dord) +1}}.
\end{align*}

\end{enumerate}
\end{theorem}


\begin{remark}
Based on Lemma \ref{lemma:L1property}, the above decay estimates extend to larger classes of functions  
 subject to the Sobolev-like constraint that $\widehat{\varphi}$ be $n_{\Dord}=\lceil \real(\Dord)\rceil +1$ times differentiable such that $\big| \widehat{\varphi}^{(j)}(\omega) \big| \big( 1+|\omega|^{j+r}\big)$ is bounded for all $0\leq j\leq n_{\Dord}$ ($r$ is a real strictly larger than $-\{-\real(\Dord)\}$).
\end{remark}

\subsection{Properties of Fractional Integration Operators}
\label{sec:propertiesintegration}

    Now that we have characterized the complex-order derivative operators, we shall look into their inverses. The difficulty there is that 
    derivative operators (integer, real or complex order) generally have null spaces which makes them non-invertible. Accordingly, our goal is to identify the class of
    integration operators that are the right inverse of a complex-order derivative as well as scale-invariant.

\begin{prop}\label{prop:Integration_Adjoint}
Let $\Dord\in \C$ with $\real(\Dord)>0$, and $a,b$ be arbitrary non-zero complex numbers. Define the integration operator $\Iop_{a,b}^{\Dord; \scalebox{0.6}{$k$}}$ of order $\Dord$ associated with the homogeneous generalized function $h_{a,b}^{\Dord}$ as in \eqref{eq:fracI_Def}, where $k\geq \max\big(1\,,\,\lfloor \real(\Dord)\rfloor \big)$. Then, 
\begin{enumerate}[(i)] 

\item $\Iop_{a,b}^{\Dord; \scalebox{0.6}{$k$}}$ defines a scale-invariant operator of degree $(-\Dord)$ over $\mathcal{S}$,

\item \label{item:left_inverse}$\Iop_{a,b}^{\Dord; \scalebox{0.6}{$k$}}$ is a right-inverse of $\Dop_{a,b}^{\Dord}$ if $k=\lfloor \real(\Dord)\rfloor$ or $k=\lceil \real(\Dord)\rceil$,

\item the adjoint operator of $\Iop_{a,b}^{\Dord; \scalebox{0.6}{$k$}}$ denoted by $\Iop_{a,b}^{(\Dord; \scalebox{0.6}{$k$})*}$ is given in \eqref{eq:fracIAdj_Def}.

\end{enumerate}
\end{prop}


\begin{remark}\label{rem:TwoAdjoints}
Since there are two admissible values of $k$ in (\ref{item:left_inverse}) for $\real(\Dord)\geq 1$ and $\real(\Dord)\not\in\N$, 
Proposition \ref{prop:Integration_Adjoint} introduces two right-inverse operators for each LSI derivative operator $\Dop_{a,b}^{\Dord}$. 
For $\real(\Dord)\in\N$ 
the two integration  operators coincide. 
The case   $0<\real (\Dord)< 1$ is  more elaborate. Let $h_{a,b}^{\Dord}(\omega)$ denote the homogeneous Fourier multiplier (frequency response) of the derivative operator $\Dop_{a,b}^{\Dord}$ with $0<\real(\Dord)< 1$. Theorem \ref{prop:Integration_Adjoint} introduces only the  integration  operator $\Iop_{a,b}^{\Dord; \scalebox{0.6}{$1$}}$ such that 
\begin{align*}
  \big( \Iop_{a,b}^{\Dord; \scalebox{0.6}{$1$}} \, \varphi \big)(x) 
= \tfrac{1}{2\pi}\int_{\R} \widehat{\varphi}(\omega) 
\tfrac{\ee^{\jj \omega x}-1 }{ h_{a,b}^{\Dord}(\omega) } \dd \omega,
\end{align*}
which is not shift-invariant (similar to all the other operators introduced in Proposition \ref{prop:Integration_Adjoint}). We should, nevertheless, mention that Theorem \ref{theo:decayRate} covers the range $\real(\Dord) > -1$. Particularly, for $-1 < \real (\Dord) < 0$, we obtain a shift-invariant integration operator. Therefore, the second right-inverse operator for $\Dop_{a,b}^{\Dord}$ when $0<\real(\Dord)< 1$ is the shift-invariant operator $\Dop_{a',b'}^{-\Dord}$:
\begin{align*}
\varphi(x) \stackrel{ \Dop_{a',b'}^{-\Dord} }{\longmapsto} \big(\Dop_{a',b'}^{-\Dord} \, \varphi \big) (x) 
= \tfrac{1}{2\pi}\int_{\R} \widehat{\varphi}(\omega) h_{a',b'}^{-\Dord}(\omega) \ee^{\jj \omega x}    \dd \omega
= \tfrac{1}{2\pi}\int_{\R} \tfrac{ \widehat{\varphi}(\omega) \ee^{\jj \omega x} }{  h_{a,b}^{\Dord}(\omega) }   \dd \omega,
\end{align*}
where $a'=\frac{1}{a}$ and $b'=\frac{1}{b}$.
\end{remark}


\begin{theorem}\label{theo:Integrator1}
Let $\Iop_{a,b}^{(\Dord; \scalebox{0.6}{$k$})*}$ be the scale-invariant operator associated with the homogeneous generalized function $h_{a,b}^{\Dord}(\omega)$ of degree $\Dord\in\C$ ($\real(\Dord)>0$) as defined in \eqref{eq:fracIAdj_Def}, such that  $k\geq \max(1\,,\,\lfloor \real(\Dord) \rfloor)$. Then, $\big( \Iop_{a,b}^{(\Dord; \scalebox{0.6}{$k$})*}\, \varphi\big)(x)$ is a well-defined and continuous function at $x\neq 0$ for all $\varphi\in\SS$. Moreover,
\begin{enumerate}[(i)]
\item \label{stat:Integrator1_NonInteger}
if $\real(\Dord)\not\in \N$ and $k=\lfloor \real(\Dord) \rfloor$, then, $\big(\Iop_{a,b}^{(\Dord; \scalebox{0.6}{$k$})*} \,\varphi\big)(x)$ is bounded around $x=0$ and the function belongs to $L_{p}$ for $p>\frac{1}{k+1-\real(\Dord)}$. More precisely, there exists a constant $c\in\R^+$ such that
\begin{align}\label{eq:IntegralDecayGood}
\forall \, x\in \R,~1\leq |x|:~~~ \Big| \big(\Iop_{a,b}^{(\Dord; \scalebox{0.6}{$k$})*} \,\varphi\big)(x)\Big| \leq \frac{c}{ |x|^{k+1-\real(\Dord)}}.
\end{align}

\item \label{stat:Integrator1_PureInt}
if $\Dord\in \N$, $h_{a,b}^{\Dord}(\omega)\equiv d\,\omega^{\Dord}$ for some $d\in\C$  and $k=\lfloor \real(\Dord) \rfloor$, then,  $\Iop_{a,b}^{(\Dord; \scalebox{0.6}{$k$})*}\, \varphi$ is bounded around $x=0$ and belongs to $L_{p}$ for all $p>0$.

\item \label{stat:Integrator1_Integer}
if $\real(\Dord)\in \N$ but $h_{a,b}^{\Dord}(\omega)\not\equiv c\,\omega^{\Dord}$  and $k=\lfloor \real(\Dord) \rfloor$, then,  $\big(\Iop_{a,b}^{(\Dord; \scalebox{0.6}{$k$})*}\, \varphi\big)(x)$ is singular at $x=0$ but $\frac{|(\Iop_{a,b}^{\Dord; \scalebox{0.4}{$k$}*}\, \varphi)(x)|}{\log|x|}$ is bounded around $x=0$. Further,  $I_{a,b}^{(\Dord; \scalebox{0.6}{$k$})*}\, \varphi$  belongs to $L_{p}$ for $p>1$, and  there exists a constant $c\in\R^+$ such that
\begin{align}\label{eq:IntegralDecaySlow}
\forall \, x\in \R,~1\leq |x|:~~~ \Big| \big(\Iop_{a,b}^{(\Dord; \scalebox{0.6}{$k$})*}\, \varphi\big)(x)\Big| \leq c\frac{1+\log|x|}{ |x|}.
\end{align}

\item \label{stat:Integrator2_NonInteger}
if $\real(\Dord)\not\in \N$ and $k>\lfloor \real(\Dord) \rfloor$, then, $|x|^{k-\real(\Dord)} |(\Iop_{a,b}^{\Dord; \scalebox{0.4}{$k$}*}\, \varphi)(x)|$ is bounded around $x=0$ (possible singularity of $\big(\Iop_{a,b}^{(\Dord; \scalebox{0.6}{$k$})*}\, \varphi\big)(x)$ at $x=0$). In addition, $\Iop_{a,b}^{(\Dord; \scalebox{0.6}{$k$})*}\, \varphi$ belongs to $L_p$ for $\frac{1}{k+1-\real(\Dord)}<p<\frac{1}{k-\real(\Dord)}$, and there exists a constant $c\in\R^+$ such that
\begin{align}\label{eq:Integral2Decay1}
\forall \, x\in \R,~1\leq |x|:~~~ \Big| \big(\Iop_{a,b}^{(\Dord; \scalebox{0.6}{$k$})*}\, \varphi\big)(x)\Big| \leq \frac{c}{ |x|^{k+1-\real(\Dord)}}.
\end{align}

\item \label{stat:Integrator2_Integer}
if $\real(\Dord)\in \N$ and $k>\lfloor \real(\Dord) \rfloor$, then, $|x|^{k-\real(\Dord)} |(\Iop_{a,b}^{\Dord; \scalebox{0.4}{$k$}*}\, \varphi)(x)|$ is bounded around $x=0$ (possible singularity of $\big(\Iop_{a,b}^{(\Dord; \scalebox{0.6}{$k$})*}\, \varphi\big)(x)$ at $x=0$). In addition, $\Iop_{a,b}^{(\Dord; \scalebox{0.6}{$k$})*}\, \varphi$ belongs to $L_p$ for $\frac{1}{k+1-\real(\Dord)}<p<\frac{1}{k-\real(\Dord)}$, and there exists a constant $c\in\R^+$ such that
\begin{align}\label{eq:Integral2Decay2}
\forall \, x\in \R,~1\leq |x|:~~~ \Big| \big(\Iop_{a,b}^{(\Dord; \scalebox{0.6}{$k$})*}\, \varphi\big)(x)\Big| \leq c\frac{1+\log|x|}{ |x|^{k+1-\real(\Dord)}}.
\end{align}

\end{enumerate}
\end{theorem}


\begin{remark}\label{rmk1}
Theorems \ref{theo:Integrator1} indicates that for $\Dord\in\C$ with $\real(\Dord)\geq 1$, it is possible to set $k$ such that $\Iop_{a,b}^{(\Dord; \scalebox{0.6}{$k$})*} : \SS \rightarrow L_p$ for any given $p>0$, except when $\frac{1}{p}+\real(\Dord)\in\Z$. For $\real(\Dord)< 1$ (even negative $\real(\Dord)$), as $k\geq 1$ is required in Theorems \ref{theo:Integrator1}, the $p$ values above $\frac{1}{1-\real(\Dord)}$ are excluded. In some sense, $k=0$ is required to cover the (almost) full range of $p>0$ values. Indeed, $k=0$ can be interpreted as
\begin{align*}
\big(\Iop_{a,b}^{\Dord; \scalebox{0.6}{$0$}}\,\varphi\big)(x) &= \mathcal{F}^{-1}\bigg\{ \frac{\widehat{\varphi}(\omega)}{ h_{a,b}^{\Dord}(\omega)}\bigg\}(x),\nonumber\\
\big(\Iop_{a,b}^{(\Dord; \scalebox{0.6}{$0$})*}\,\varphi\big)(x) &= \mathcal{F}^{-1}\bigg\{ \frac{\widehat{\varphi}(\omega)}{ h_{a,b}^{\Dord}(-\omega)}\bigg\}(x).
\end{align*}
This suggests that $\Iop_{a,b}^{\Dord; \scalebox{0.6}{$0$}} = \Dop_{a',b'}^{-\Dord}$ and $\Iop_{a,b}^{(\Dord; \scalebox{0.6}{$0$})*} = \Dop_{b',a'}^{-\Dord}$, where $\real(-\Dord)>-1$, and $a'=\frac{1}{a},\,b'=\frac{1}{b}$. Theorem \ref{theo:decayRate} shows that $\Iop_{a,b}^{(\Dord; \scalebox{0.6}{$0$})*}\,\varphi$ is in $L_p$ for $p<\frac{1}{1-\real(\Dord)}$. Thus, by extending the family of $\Dord$-order integration operators and adjoints to include $k=0$ as above (for $\real(\Dord)< 1$), we are always able to choose a $k$ such that $\Iop_{a,b}^{(\Dord; \scalebox{0.6}{$k$})*} : \SS \rightarrow L_p$, except when $\frac{1}{p}+\real(\Dord)\in\Z$.
\end{remark}


\subsection{Impulse Responses of Fractional Derivative and Integration Operators}
To  evaluate the output of complex-order operators, we start by simple inputs such as Dirac's delta function and its shifted versions. This then yields the Schwartz kernel of these operators.

\begin{theorem}\label{theo:ImpulseResponse}
Let $\Dop_{a,b}^{\Dord}$, $\Iop_{a,b}^{\Dord; \scalebox{0.6}{$k$}}$ and $\Iop_{a,b}^{(\Dord; \scalebox{0.6}{$k$})*}$ be the scale-invariant operators as defined in \eqref{eq:fracD_Def}, \eqref{eq:fracI_Def} and \eqref{eq:fracIAdj_Def}, respectively, where $k\geq \max(1\,,\,\lfloor \real(\Dord)\rfloor)$ and $\Dord\not\in\Z$.
Then, the response of $\Dop_{a,b}^{\Dord}$, $\Iop_{a,b}^{\Dord; \scalebox{0.6}{$k$}}$ and $\Iop_{a,b}^{(\Dord; \scalebox{0.6}{$k$})*}$ to the shifted Dirac impulse $\delta(\cdot-\tau)$ is given by 
\begin{align}\label{eq:ImpulseResponse}
\Dop_{a,b}^{\Dord}\big\{ \delta(\cdot-\tau)\big\}(x)
=\phantom{-}& 
\tfrac{\Gamma(\Dord+1)}{2\pi}\Big(\tfrac{a}{(\jj \tau -\jj x)^{\Dord+1}} + \tfrac{b}{(\jj x-\jj \tau)^{\Dord+1}} \Big), \nonumber\\
\Iop_{a,b}^{\Dord; \scalebox{0.6}{$k$}} \big\{\delta(\cdot-\tau)\big\}(x) 
= \phantom{-}& 
\tfrac{\Gamma(1-\Dord)}{2\pi}\Big( \tfrac{(\jj \tau-\jj x)^{\Dord-1}}{a} + \tfrac{(\jj x-\jj \tau)^{\Dord-1}}{b} \Big) 
\nonumber\\
-& \tfrac{\Gamma(1-\Dord)}{2\pi} \Big( \tfrac{(\jj \tau)^{\Dord-1}}{a} + \tfrac{(-\jj \tau)^{\Dord-1}}{b}\Big)\sum_{j=0}^{k-1} \scalebox{1}{$\binom{j-\Dord}{j}$} \big(\tfrac{x}{\tau}\big)^{j} , \nonumber\\
\Iop_{a,b}^{(\Dord; \scalebox{0.6}{$k$})*} \big\{\delta(\cdot-\tau)\big\}(x) 
= \phantom{-}& 
\tfrac{\Gamma(1-\Dord)}{2\pi}\Big( \tfrac{(\jj x-\jj \tau)^{\Dord-1}}{a}  + \tfrac{(\jj \tau-\jj x)^{\Dord-1}}{b} \Big) 
\nonumber\\
-& \tfrac{\Gamma(1-\Dord)}{2\pi} \Big( \tfrac{(\jj x)^{\Dord-1}}{a} + \tfrac{(-\jj x)^{\Dord-1}}{b} \Big)\sum_{j=0}^{k-1} \scalebox{1}{$\binom{j-\Dord}{j}$} \big(\tfrac{\tau}{x}\big)^{j} ,
\end{align}
where $\binom{x}{0}=1$ and $\binom{x}{j} = \tfrac{\Gamma(x+1)}{j! \Gamma(x+1-j)} = \tfrac{x(x-1)\dots(x-j+1)}{j!}$ extends the standard definition of $\binom{n}{j}$.
\end{theorem}


\begin{remark}
Theorem \ref{theo:ImpulseResponse} enables us to represent the operators in terms of integrals. Indeed, for $\varphi\in\mathcal{S}$, we know that
\begin{align*}
\varphi(x) = \int_{\R} \varphi(\tau)\,\delta(x-\tau) \dd\tau.
\end{align*} 
According to  Schwartz' kernel theorem \cite{Schwartz1950,Treves1967}, the action of $\Lop:\mathcal{S} \to \mathcal{S}'$ can be represented as 
\begin{align*}
\big(\Lop \varphi\big)(x) = \int_{\R} \varphi(\tau)\,\big(\Lop\delta(\cdot-\tau)\big)(x) \dd\tau, \nonumber\\
\end{align*} 
where $\Lop$ stands for any of $\Dop_{a,b}^{\Dord}$, $\Iop_{a,b}^{\Dord; \scalebox{0.6}{$k$}}$ and $\Iop_{a,b}^{(\Dord; \scalebox{0.6}{$k$})*}$. The integral forms are particularly useful for numerical evaluation of $\big(\Lop \varphi\big)(x)$.
\end{remark}


\begin{remark}
With the particular choice of $a= \jj^{\Dord-1}$ and $b=(-\jj)^{1-\Dord}$, the kernels in Theorem \ref{theo:ImpulseResponse} simplify to more familiar forms of 
\begin{align}\label{eq:ImpulseResponseSimplified}
\big(\Dop_{a,b}^{\Dord} \delta(\cdot-\tau)\big)(x) =\phantom{-}& 
\tfrac{\Gamma(1+\Dord)\sin(\pi\Dord)}{\pi} (x-\tau)_{+}^{-1-\Dord}, \nonumber\\
\big(\Iop_{a,b}^{\Dord; \scalebox{0.6}{$k$}} \delta(\cdot-\tau)\big)(x) = \phantom{-}& 
\tfrac{\Gamma(1-\Dord)\sin(\pi\Dord)}{\pi} \Big( (x-\tau)_{+}^{\Dord-1} - (-\tau)_{+}^{\Dord-1}  \sum_{j=0}^{k-1} \scalebox{1}{$\binom{j-\Dord}{j}$} \big(\tfrac{x}{\tau}\big)^{j} \Big),
\nonumber\\
\big(\Iop_{a,b}^{(\Dord; \scalebox{0.6}{$k$})*} \delta(\cdot-\tau)\big)(x) = \phantom{-}& 
\tfrac{\Gamma(1-\Dord)\sin(\pi\Dord)}{\pi} \Big( (\tau-x)_{+}^{\Dord-1} - (-x)_{+}^{\Dord-1}  \sum_{j=0}^{k-1} \scalebox{1}{$\binom{j-\Dord}{j}$} \big(\tfrac{\tau}{x}\big)^{j} \Big).
\end{align}
\end{remark}


\section{Self-Similar Stable Processes with Complex Hurst Index}
\label{sec:processes}

    Our goal in this section is to specify a class of self-similar stable processes with complex-valued Hurst exponent. 
    The random processes as specified as generalized random processes via their characteristic functional, both concepts being reintroduced in Section~\ref{sec:GRPandCF}. 
    We then use the integration operators defined in Section~\ref{sec:Theos} to construct our family of self-similar stable processes in Section~\ref{sec:constructprocess}. 
    The invariance and regularity properties of self-similar stable processes are studied in Sections~\ref{sec:propertiesofprocesses} and \ref{sec:regularityofstable} respectively. Finally, we provide simulations of realizations of the considered class of random processes in Section~\ref{sec:simul}. 

    \subsection{Generalized Random Processes and their Characteristic Functional} \label{sec:GRPandCF}
    
     The theory of generalized random processes has been formalized by Gelfand and It\^o~\cite{Gelfand_Vilenkin,Ito1984foundations}.
    It is the probabilistic counterpart of the theory of generalized functions of Schwartz~\cite{Schwartz1966distributions} and allows for the construction of broad classes of random processes, including the ones that do not admit pointwise representations such as the stable white noises~\cite[Chapter 3]{Gelfand_Vilenkin}. We briefly recap this formalism, with a special emphasis on the characteristic functional. 
    
A tempered generalized random process $S$ (or simply a generalized random process) is a collection of random variables $\langle S , \varphi \rangle$, where $\varphi$ is a test function in the Schwartz space $\mathcal{S}$, such that 
\begin{itemize}
    \item \textit{Linearity:} for any $\varphi_1, \varphi_2 \in \mathcal{S}$ and $\lambda \in \R$, $\langle S , \varphi_1 +  \lambda \varphi_2\rangle = \langle S , \varphi_1 \rangle + \lambda \langle S , \varphi_2 \rangle$ almost surely; 
    
    \item \textit{Continuity:} for any converging sequence $\varphi_n \rightarrow \varphi$ in $\mathcal{S}$, the random variables $\langle S  , \varphi_n \rangle$ converge in probability to $\langle S , \varphi \rangle$. 
\end{itemize}

The space $\mathcal{S}'$ of tempered generalized functions is endowed with the weak* topology \cite{Treves1967}. This topology defines a Borelian $\sigma$-field on $\mathcal{S}'$ and a tempered generalized random process $S$ can then be seen as a random element of the space of tempered generalized functions. The probability law of a generalized random process $S$ is the probability measure $\mathscr{P}_S$ over $\mathcal{S}'$ such that, for any Borelian set $B \subset \mathcal{S}'$, we have $\mathscr{P}_S (B) = \mathscr{P} (S\in B)$. 
We refer the interested reader to~\cite{Bierme2017generalized} for a comprehensive introduction to these notions in the framework of tempered generalized functions and for additional references. 

We are interested in complex-valued random processes; hence, we shall adapt the usual concepts to this case. Thereafter, $\mathcal{S}$ and $\mathcal{S}'$ denote the complex-valued Schwartz space and space of tempered generalized functions, respectively. The characteristic function of a complex random variable $Z$ is given by $\widehat{\mathscr{P}}_Z(\xi) = \mathbb{E} [ \mathrm{e}^{\mathrm{i} \real( Z \bar{\xi} ) }]$ where $\bar{\xi}$ is the complex conjugate of $\xi\in \C$. The extension to complex-valued generalized random processes is as follows. 

\begin{definition}[\bf Characteristic functional]
The characteristic functional of a complex valued tempered generalized random process $S$ is defined as 
\begin{equation} \label{eq:CFS}
    \widehat{\mathscr{P}}_S (\varphi) = \mathbb{E} \left[ \mathrm{e}^{\mathrm{i} \real( \langle S , \bar{\varphi} \rangle) } \right] = \int_{\mathcal{S}'} \mathrm{e}^{ \mathrm{i} \real( \langle u , \bar{\varphi} \rangle ) } \mathrm{d}\mathscr{P}_S(u)
\end{equation}
for any $\varphi \in \mathcal{S}$.
\end{definition}

The characteristic functional has been introduced by Kolmogorov~\cite{Kolmogorov1935transformation} and popularized by Gelfand and Vilenkin~\cite{Gelfand_Vilenkin}.
The probability law of a generalized random process is characterized by its characteristic functional, which is its (infinite-dimensional) Fourier transform. 
In particular, two random processes with identical characteristic functionals have identical finite-dimensional marginals. 

The Bochner-Minlos theorem ensures that the characteristic functionals of tempered generalized processes are exactly the functionals $\widehat{\mathscr{P}} : \mathcal{S}\rightarrow \C$ that are continuous, positive-definite, and such that $\widehat{\mathscr{P}}(0) = 1$~\cite{Gelfand_Vilenkin}. This provides a way of constructing tempered generalized processes via their characteristic functionals. 

For example, for any $\alpha \in ]0,2]$, it is known that the functional 
\begin{equation}
    \widehat{\mathscr{P}} (\varphi) = 
    \exp \left( - \int_{\R} |\varphi( x )|^\alpha \mathrm{d}x \right) =
    \exp \left( - \lVert \varphi \rVert_\alpha^\alpha \right)
\end{equation}
is continuous and positive-definite over real-valued Schwartz functions and satisfies $\widehat{\mathscr{P}} (\varphi) = 0$~\cite{Fageot2014}. Therefore, there exists a real generalized random process $W_\alpha$, called a symmetric-$\alpha$-stable (S$\alpha$S) white noise, such that $\widehat{\mathscr{P}}_{W_\alpha} (\varphi) = \widehat{\mathscr{P}} (\varphi)$. 
The family of stable white noises $W_\alpha$ will play an important role in this paper. Note that, seen as a complex-valued generalized random process acting on complex-valued test functions $\varphi \in \mathcal{S}$ and using that $W_\alpha$ itself is real, we have, due to \eqref{eq:CFS}, that
\begin{equation}
\widehat{\mathscr{P}}_{W_\alpha} (\varphi)  = \mathbb{E} \left[ \mathrm{e}^{\mathrm{i} \real( \langle W_\alpha , \bar{\varphi} \rangle)} \right] = \exp( - \lVert \real( \varphi ) \rVert_\alpha^\alpha).    
\end{equation}

    \subsection{Construction of Complex-valued Self-Similar Stable Processes}
    \label{sec:constructprocess}
    
Let $\alpha \in ]0,2]$. 
We can construct new classes of generalized random processes as follows. 
Assume that $\mathrm{U}$ is a linear and continuous operator from $\mathcal{S}$ to $L_\alpha$, where both spaces are for complex valued (generalized) functions\footnote{The space $L_\alpha = L_{\alpha}( \R , \C)$ is the space of measurable functions $f : \R \rightarrow \C$ whose real and imaginary parts are in $L_\alpha(\R,\R)$.}. 
The properties of both $W_\alpha$ and $\mathrm{U}$ ensure that $\varphi \mapsto \exp\left( - \lVert \real( \mathrm{U} \{\varphi \} ) \rVert_\alpha^\alpha \right)$ is continuous, positive-definite over $\mathcal{S}$ and normalized.
Hence, due to the Bochner-Minlos theorem, there exists a generalized random process $S$ with characteristic functional $\widehat{\mathscr{P}}_S(\varphi) = \widehat{\mathscr{P}}_{W_\alpha}( \mathrm{U} \{ \varphi \}) = \exp\left( - \lVert \real( \mathrm{U} \{ \varphi \} ) \rVert_\alpha^\alpha \right)$. We use this principle to construct our extended class of generalized random processes.

\begin{prop} \label{prop:constructS}
    Let $a,b \in \C \backslash \{0\}$, $\gamma \in \C$ with $\real(\gamma) > 0$, and $\alpha \in ]0,2]$. 
    We further assume that 
    \begin{equation} \label{eq:notinN}
        \tfrac{1}{\alpha} + \real(\gamma) \notin \N
    \end{equation} 
    and we set 
    \begin{equation} \label{eq:kalgam}
    k(\alpha,\gamma)= \left\lfloor \tfrac{1}{\alpha} + \real(\gamma) \right\rfloor \in \N. 
\end{equation}
    Then, the linear operator $(\Iop^{\Dord; \scalebox{0.6}{$k$}(\alpha,\gamma)}_{a,b})^*$ is continuous from $\mathcal{S}$ to $L_\alpha$ and there exists a tempered generalized random process $S_{a,b}^{\gamma,\alpha}$ such that
    \begin{equation}
        \widehat{\mathscr{P}}_{S_{a,b}^{\gamma,\alpha}}(\varphi) =  
        \exp\left( - \left\lVert \real\left(  (\Iop^{\Dord; \scalebox{0.6}{$k$}(\alpha,\gamma)}_{a,b})^*  \{ \varphi \}  \right) \right\rVert_\alpha^\alpha \right).
    \end{equation}
    If moreover the condition
    \begin{equation} \label{eq:conditionalpha}
        \alpha > \left\{\begin{array}{cl}
            \frac{1}{2 - \{\real(\Dord)\}} & \text{if } \real(\Dord) \not\in \mathbb{N}, \phantom{\Big|}\\
            1  & \text{if } \real(\Dord) \in \mathbb{N}, \phantom{\Big|}
        \end{array}\right.
    \end{equation}
    is satisfied, then, the random process $S_{a,b}^{\gamma,\alpha}$ can be whitened in the sense that $\Dop_{a,b}^\gamma S_{a,b}^{\gamma,\alpha} = W_\alpha$ is a S$\alpha$S white noise.
\end{prop}

\begin{proof}
    The continuity of $(\Iop^{\Dord; \scalebox{0.6}{$k$}(\alpha,\gamma)}_{a,b})^*$ for $k = k(\alpha,\gamma)\geq 1$ is a direct consequence of Theorem~\ref{theo:Integrator1}, for which condition~\eqref{eq:notinN} is required. 
    For the case of $k(\alpha,\gamma) =0$,
    which necessitates  $\real(\gamma)< 1$, we recall that $(\Iop_{a,b}^{\gamma;0})^*$ is a true (left- and right-) inverse for $\Dop_{a,b}^\gamma$. This case is not covered by Theorem~\ref{theo:Integrator1} but is discussed in Remark~\ref{rmk1}. 
    Finally, 
    the existence result simply follows from the application  of the Bochner-Minlos theorem to $\varphi \mapsto \exp\left( - \lVert \real(  (\Iop^{\Dord; \scalebox{0.6}{$k$}(\alpha,\gamma)}_{a,b})^* \varphi ) \rVert_\alpha^\alpha \right)$.
    
    For the claim that $S_{a,b}^{\gamma,\alpha}$ can be whitened using $\Dop_{a,b}^\gamma$, it is sufficient that $(\Iop_{a,b}^{\gamma;k})$ is a right-inverse of $\Dop_{a,b}^\gamma$. By Proposition \ref{prop:Integration_Adjoint}-(\ref{item:left_inverse}),  this property holds if  $k=\lfloor \real(\Dord)\rfloor$ or $k=\lceil \real(\Dord) \rceil$. Since we have $k=\lfloor \frac{1}{\alpha} + \real(\Dord)\rfloor$, this implies that
    \begin{align}\label{eq:floor_ceil}
        \lfloor \tfrac{1}{\alpha} + \real(\Dord)\rfloor \leq \lceil \real(\Dord)\rceil. 
    \end{align}
    It is then not difficult to see that \eqref{eq:floor_ceil} is equivalent to 
    \eqref{eq:conditionalpha}.
\end{proof}
    
    We call the generalized random process $S_{a,b}^{\gamma,\alpha}$ of Proposition~\ref{prop:constructS} a \textit{(complex-valued) fractional stable process}.

\begin{figure}[tb]
		\centering
	\includegraphics[width=12cm]{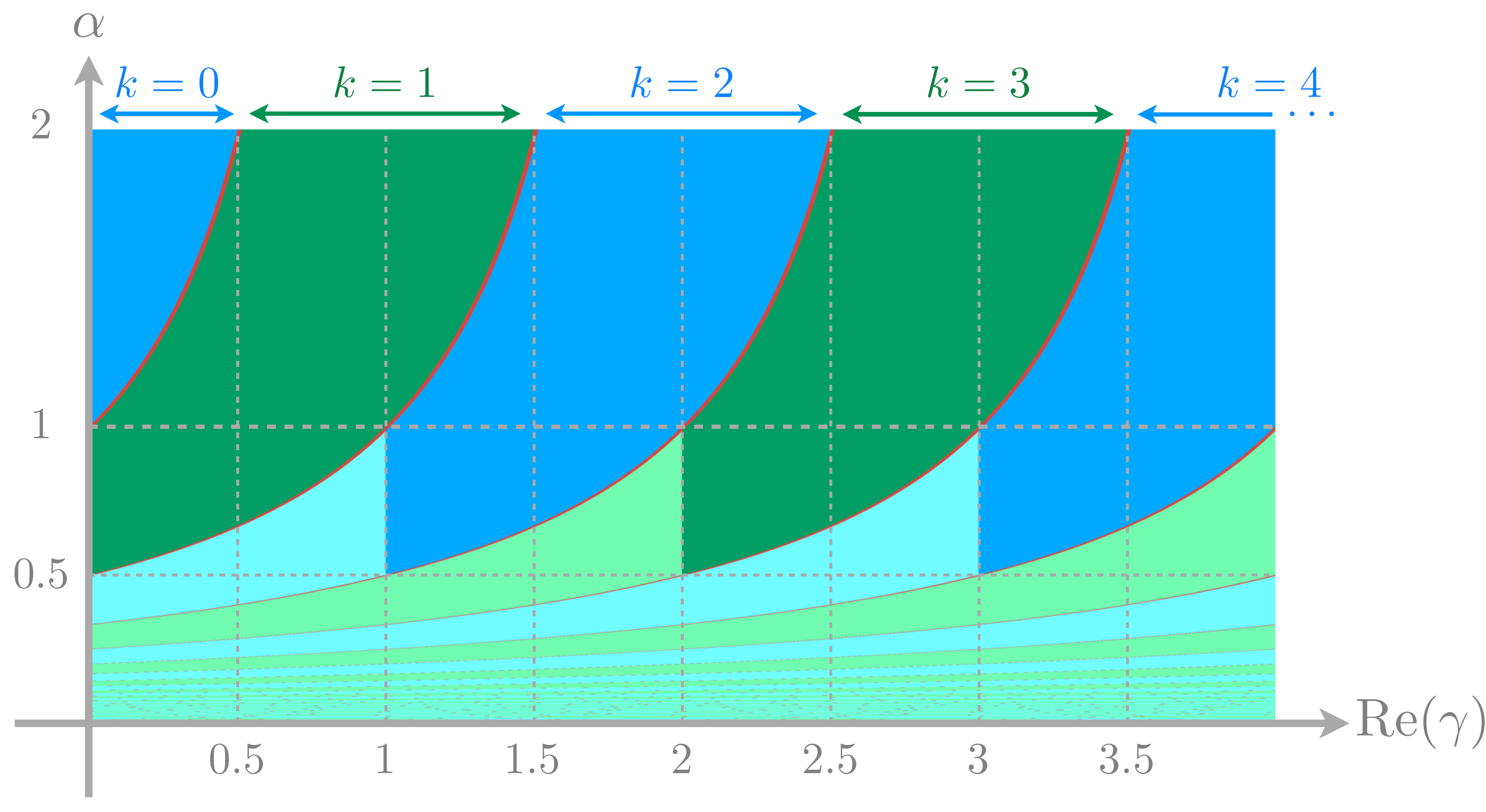}
		\caption{The  fractional stable process $S_{a,b}^{\Dord,\alpha}$ is well-defined by applying the operator $\Iop_{a,b}^{\Dord , k}$ to the symmetric-$\alpha$-stable white noise with $k= \lfloor \frac{1}{\alpha} + \real(\Dord) \rfloor$, whenever $\frac{1}{\alpha}+\real(\Dord)\not\in\mathbb{N}$ (the red curves separating the blue and green areas). For $\alpha > \frac{1}{2 - \{\real(\Dord)\}}$ (dark blue and dark green areas), the fractional process can be whitened by applying the fractional derivative operator  $\Dop_{a,b}^{\Dord}$. The use of green and blue colors is to better highlight the red curves separating them, and do not encode any other information.} \label{fig:alpha_gamma}
	\end{figure}

Figure~\ref{fig:alpha_gamma} delineates the type of  parameter pairs of $(\alpha,\gamma)$ for which the fractional stable random processes $S_{a,b}^{\gamma,\alpha}$ is well-defined and can be whitened. Note that  $\Dop_{a,b}^\gamma S_{a,b}^{\gamma,\alpha}$ is not a white process when \eqref{eq:conditionalpha} is not satisfied. However, fractional random processes are always whitenable for $\alpha > 1$. 
    
\begin{remark}
It is worth noting that the considered class of fractional stable processes are built over the family of symmetric alpha-stable S$\alpha$S white noise. Consequently, they are symmetric as well ($S_{a,b}^{\gamma,\alpha}$ and $-S_{a,b}^{\gamma,\alpha}$ have the same law). It is possible to extend the family by considering  non-symmetric white noises, which are fully described by four parameters~\cite[Definition 1.1.6]{SamorodnitskyBook94}. This extension is readily achievable for $\alpha \neq 1$, but is more technical for specific cases corresponding to $\alpha = 1$, for which the characteristic exponent $\Psi$ can include a logarithmic term. Due to the technicalities, we do not delve into the details in this paper.
\end{remark}

    \subsection{Invariance Properties of Complex-Valued Fractional Stable Processes}
    \label{sec:propertiesofprocesses}

We first adapt the notion of self-similarity to complex-valued generalized random processes that may have no pointwise interpretation.  

\begin{definition}[\bf Self-similar generalized process]
A complex-valued tempered generalized random process $S$ is self-similar with Hurst index $H \in \C$ if the probability laws of $S$ and $T^{-H}S (T \cdot)$ are identical for any $T > 0$. 
\end{definition}

When $S$ admits a pointwise interpretation, we readily recover the conventional self-similarity defined in \eqref{eq:SSprocDef}. In the present setting where the random process $S$ is complex,  the Hurst index is generally complex as well. 
However, the driving  S$\alpha$S white noise $W_\alpha$ is self-similar and real-valued with (real) Hurst index $H = \frac{1}{\alpha} - 1$ \cite[Proposition 4.2]{Fageot2019scaling}.


\begin{prop}
    \label{prop:scaleinv}
    Under the   conditions  of Proposition~\ref{prop:constructS}, the fractional stable process $S_{a,b}^{\gamma,\alpha}$ is self-similar with Hurst exponent 
    \begin{equation} \label{eq:Hselfsim}
        H= \gamma + \frac{1}{\alpha} - 1 \in \C.
    \end{equation}
\end{prop}    

\begin{proof}
    To simplify the notations, we set $S = S_{a,b}^{\gamma,\alpha}$ and $\Iop = \Iop^{\Dord; \scalebox{0.6}{$k$}(\gamma,\alpha)}_{a,b}$.
    The self-similarity can be proved using the characteristic functional and the fact that $\widehat{\mathscr{P}}_{S_1} = \widehat{\mathscr{P}}_{S_2}$ if and only if $S_1$ and $S_2$ have the same probability law. Indeed, for $T > 0$, we have that
    \begin{align}
        \widehat{\mathscr{P}}_{S(T\cdot)} (\varphi) &= \widehat{\mathscr{P}}_S\big( T^{-1} \varphi ( \cdot / T)  \big)
        = \mathrm{e}^{- \lVert T^{-1} \real(\Iop^* \{ \varphi(\cdot /T) \} ) \rVert_\alpha^\alpha } \\
        &= \mathrm{e}^{- \lVert T^{-1-\gamma} \real(\Iop^* \{ \varphi \} (\cdot /T) ) \rVert_\alpha^\alpha }
       = \mathrm{e}^{- \lVert T^{-1-\gamma + \frac{1}{\alpha}} \real(\Iop^* \{ \varphi \} ) \rVert_\alpha^\alpha }  \\
       &= \widehat{\mathscr{P}}_S(T^{-1-\gamma + \frac{1}{\alpha}} \varphi )
       =  \widehat{\mathscr{P}}_{T^{-1-\gamma + \frac{1}{\alpha}} S} (\varphi),
     \end{align}
     where we used in particular the $(-\gamma)$-homogeneity of $\Iop^*$ for the third equality. 
     Hence, $T^{-H} S(T\cdot) = S$ for any $T> 0$ where $H$ is given by \eqref{eq:Hselfsim}. 
\end{proof}

    We now study the invariance of fractional stable processes with respect to shift operations. Again, dealing \textit{a priori} with processes with no pointwise interpretation, we shall adapt the usual notions to generalized random processes. For $h_0 > 0$, we define the operator $\Delta_{h_0} \{f\} = f(\cdot + h_0) - f$, which represents the increments of a generalized function $f \in \mathcal{S}'$. 
    
    \begin{definition} \label{def:statio}
    A generalized random process $S$ is stationary if $s(\cdot - x_0)$ and $s$ have the same law for any $x_0 \in \R$. We say moreover that $S$ has stationary increments of order $k \geq 1$ if the generalized random processes $(\Delta_{h_0})^k S$ are stationary for any $h_0>0$.
    \end{definition}
    
    \begin{remark}
    Classically, we say that a random process $S = (S(x))_{x \in \R}$ with well-defined sampled values has stationary increments if the law of $s(x_1) - s(x_0)$ only depends on the difference $(x_1 - x_0)$ for any $x_0 < x_1$. Definition~\ref{def:statio} covers and generalizes this notion. Indeed, assume that $s$ has stationary increments of order $k=1$ in the sense of Definition~\ref{def:statio}. Then, for any test function $\varphi$, $h_0 > 0$ and $t \in \R$, we have that $\langle \Delta_{h_0}S , \varphi \rangle$ and $\langle \Delta_{h_0} S , \varphi (\cdot - t) \rangle$ are equal in law. Picking $\varphi = \delta(\cdot - t_0)$ (which is possible for pointwise random processes), $h_0 = t_1 - t_0$, and $t = -t_0$, we deduce that
    \begin{equation}
        s(t_1) - s(t_0) \overset{(\mathcal{L})}{=} s(t_1-t_0) - s(0).
    \end{equation}
    The latter only depends on $t_1 - t_0$ and $S$ has stationary increments in the classical sense.
    \end{remark}
    
    \begin{prop}
 Under the   conditions of  Proposition~\ref{prop:constructS}, the fractional stable process $S_{a,b}^{\gamma,\alpha}$ 
    \begin{itemize}
        \item is stationary if $0 < \real(\gamma) < 1$; and 
        \item has stationary increments of order $\lfloor \real(\gamma)  \rfloor$ if $\real(\gamma) \geq 1$. 
    \end{itemize}
 \end{prop}

\begin{proof}
We set $S = S_{a,b}^{\gamma,\alpha}$ and $\Iop = \Iop^{\Dord; \scalebox{0.6}{$k$}(\gamma,\alpha)}_{a,b}$.
    The proof relies on the characteristic functional. For $\real(\gamma) \in ]0,1[$, we have that $k(\gamma, \alpha) = 0$ (see \eqref{eq:kalgam}) and 
    the operator $\Iop^* = (\Iop_{a,b}^{\gamma; 0})^*  = ((\Dop_{a,b}^{\gamma})^{-1})^*$ is a convolution; hence 
    \begin{align} \label{eq:CFstatio}
        \CF_{S(\cdot - x_0)} (\varphi) 
        &= \CF_{S}(\varphi(\cdot + x_0)  )
        = \exp( - \lVert \real(\Iop^* \{ \varphi (\cdot + x_0)  ) \} \rVert_\alpha^\alpha )   \nonumber \\
        &= \exp( - \lVert \real(\Iop^* \{ \varphi \} (\cdot + x_0) )  \rVert_\alpha^\alpha ) 
        = \exp( - \lVert \real( \Iop^* \{ \varphi \} )  \rVert_\alpha^\alpha )   \nonumber \\
        &= \CF_{S} (\varphi) ,
    \end{align}
    where we used the shift-invariance of $\Iop^*$ in the third equality and a simple change of variables in the fourth one. Hence, $S$ is stationary. \\
    
    Assume that $\real(\gamma) \geq 1$. We treat the case $\real(\gamma) \in [1,2[$. First, using \eqref{eq:fracIAdj_Def} with $\Delta_{- h_0} \varphi$, we remark that $\widehat{\Delta_{- h_0} \varphi}(0) = \widehat{\varphi(\cdot - h_0)} (0) - \widehat{\varphi}(0)=  0$ and therefore
    \begin{equation}
        \Iop^*(\Delta_{-h_0}\varphi) = \mathcal{F}^{-1} \left(\frac{ \widehat{\Delta_{-h_0} \varphi} - \widehat{\Delta_{-h_0} \varphi}(0)}{h_{a,b}^\gamma} \right)
        = \mathcal{F}^{-1} \left(\frac{ \widehat{\Delta_{-h_0} \varphi}  }{h_{a,b}^\gamma} \right) = (\Dop_{a,b}^{\gamma,*})^{-1} (\Delta_{-h_0} \varphi). 
    \end{equation}
    In particular, when restricted to functions of the form $\Delta_{-h_0}\varphi$, $\Iop^*$ is a convolution. 
    As we did in \eqref{eq:CFstatio}, we  prove that 
    $\CF_{(\Delta_{h_0} S)(\cdot - x_0)} (\varphi) = \CF_{\Delta_{h_0} S} (\varphi)$ 
    for any $\varphi \in \mathcal{S}$. 
    This shows that $\Delta_{h_0} S$ is stationary, 
    or equivalently, $S$ has stationary increments of order $1 = \lfloor \gamma \rfloor$. 
    
    For $\gamma \geq 2$, we set $k= \lfloor \gamma \rfloor$. 
    Then, the function $\psi = (\Delta_{-h_0})^k \varphi$ is such that $\widehat{\psi}^{(j)} (0) = 0$ for any $0\leq j \leq k-1$. Hence, $\Iop^* \psi = (\Dop_{a,b}^{\gamma *})^{-1} \psi$ and the same argument as for $\gamma \in [1,2[$ applies.
\end{proof}
    
    \subsection{Regularity of Fractional Stable Processes}
    \label{sec:regularityofstable}
    
    We characterize the smoothness of fractional stable processes in terms of local Sobolev regularity in the space $W_{p, \mathrm{loc}}^{\tau}$ with $\tau \in \R$ and $1 \leq p \leq \infty$. With $p = 2$, we recover the $L_2$-Sobolev regularity, while $p = \infty$ corresponds to the Hölder regularity~\cite{Triebel2010theory}. For fixed $p\geq 1$ and $\tau_1 \leq \tau_2$, we have the continuous embedding $W_{p,\mathrm{loc}}^{\tau_2} \subseteq W_{p,\mathrm{loc}}^{\tau_1}$. 
    
    Following~\cite{fageot2020critical}, we characterize the regularity properties of a generalized (random) function $f \in \mathcal{S}'$ via its critical smoothness function, defined for $1 \leq p \leq \infty$ by
    \begin{equation}
        \tau_f(p) = \sup \{ \tau \in \R, \ f \in W_{p, \mathrm{loc}}^{\tau} \}.
    \end{equation}
    Then, $f \in W_{p, \mathrm{loc}}^{\tau}$ for any $\tau < \tau_f(p)$ and $f \notin W_{p, \mathrm{loc}}^{\tau}$ for  $\tau > \tau_f(p)$.

    \begin{theorem}
 Under the   conditions of  Proposition~\ref{prop:constructS}, the fractional stable process $S_{a,b}^{\gamma,\alpha}$ has the following properties.
 \begin{itemize}
     \item If $\alpha = 2$ (Gaussian case), then $\tau_{S_{a,b}^{\gamma,\alpha}}(p) =\real(\gamma)\ - \frac{1}{2}$.
     \item If $\alpha < 2$, then $\tau_{S_{a,b}^{\gamma,\alpha}}(p) = \real(\gamma) + \frac{1}{\max(p,\alpha)} - 1$.
 \end{itemize}
    \end{theorem}
    
    \begin{proof}
     The proof follows from the combination of the  two facts.
    First, the critical Sobolev smoothness $\tau_{W_\alpha}(p)$ of a S$\alpha$S white noise has been fully characterized in the Gaussian case in~\cite{Veraar2010regularity} and in the general case in a series of papers~\cite{Fageot2017besov,Fageot2017multidimensional,aziznejad2020wavelet}. For any $p\geq 1$, the results are synthesized in~\cite[Theorem 1]{aziznejad2020wavelet} as
    \begin{align}\label{eq:regunoise_2}
        \tau_{W_2}(p) = - 1 / 2
    \end{align}
    for the Gaussian ($\alpha=2$) case and
    \begin{align}\label{eq:regunoise_p}
        \tau_{W_\alpha}(p) = \frac{1}{\max(p,\alpha)} - 1
    \end{align}
    for the non-Gaussian ($\alpha < 2$) case.

    Second, we show that the operators $\Dop_{a,b}^{\Dord}$ induce a systematic decrease of the Sobolev smoothness in the sense that $f \in W_{p, \mathrm{loc}}^{\tau}$ if and only if $\Dop_{a,b}^{\Dord} f \in W_{p, \mathrm{loc}}^{\tau- \mathrm{Re}(\gamma)}$ for any $p\geq 1$ and $\tau \in \R$. This is proven by applying the criterion of~\cite[Theorem 2]{Fageot2017nterm} to this setting\footnote{We observe that~\cite{Fageot2017nterm} deals with periodic random processes. The results easily apply to our setting since we consider the \textit{local} regularity of the proposed self-similar stable processes.}. Indeed, denoting by $m(\omega) = \frac{h_{a,b}^\gamma(\omega)}{|\omega|^{\mathrm{Re}(\gamma)}}$, we observe that $|m(\omega)| = |a|$ if $\omega > 0$ and $|m(\omega)| = |b|$ for $\omega < 0$ and the relation~\cite[Eq. (28)]{Fageot2017nterm} is readily satisfied. 
    
    Then, the integration operator $\Iop^{\Dord; \scalebox{0.6}{$k$}(\gamma,\alpha)}_{a,b}$, which is a right-inverse of $\Dop_{a,b}^{\Dord}$,  satisfies the converse relation that $f \in W_{p, \mathrm{loc}}^{\tau}$ if and only if $g = \Iop^{\Dord; \scalebox{0.6}{$k$}(\gamma,\alpha)}_{a,b} f \in W_{p, \mathrm{loc}}^{\tau + \mathrm{Re}(\gamma)}$ (using that $\Dop_{a,b}^\gamma g = f$). This shows that the critical smoothness of  $S_{a,b}^{\gamma,\alpha} = \Iop^{\Dord; \scalebox{0.6}{$k$}(\gamma,\alpha)}_{a,b} W_\alpha$ is such that
    \begin{equation}
        \tau_{S_{a,b}^{\gamma,\alpha} }(p) = \tau_{W_\alpha} (p) + \mathrm{Re}(\gamma),
    \end{equation}
    and the result follows from the stable white noise case in \eqref{eq:regunoise_2} and \eqref{eq:regunoise_p}.
    \end{proof}

    \begin{remark}
    The critical $L_2$-Sobolev regularity, corresponding to $p=2$, does not depend on $\alpha$ and is equal to $ \tau_{S_{a,b}^{\gamma,\alpha}}(2) = \real(\gamma) - \frac{1}{2}$ for any fractional stable process. Moreover, the critical Hölder regularity is $\real(\gamma) - \frac{1}{2}$ for $\alpha = 2$ and $\real(\gamma) - 1$ otherwise. This last result is coherent with \cite{Huang2007} for $\gamma \geq 1$ and generalizes the result for new fractional stable random processes. 
    \end{remark}

    
    \subsection{Simulations}
        \label{sec:simul}
    
\begin{figure}[tb]
		\centering
	\includegraphics[width=12cm]{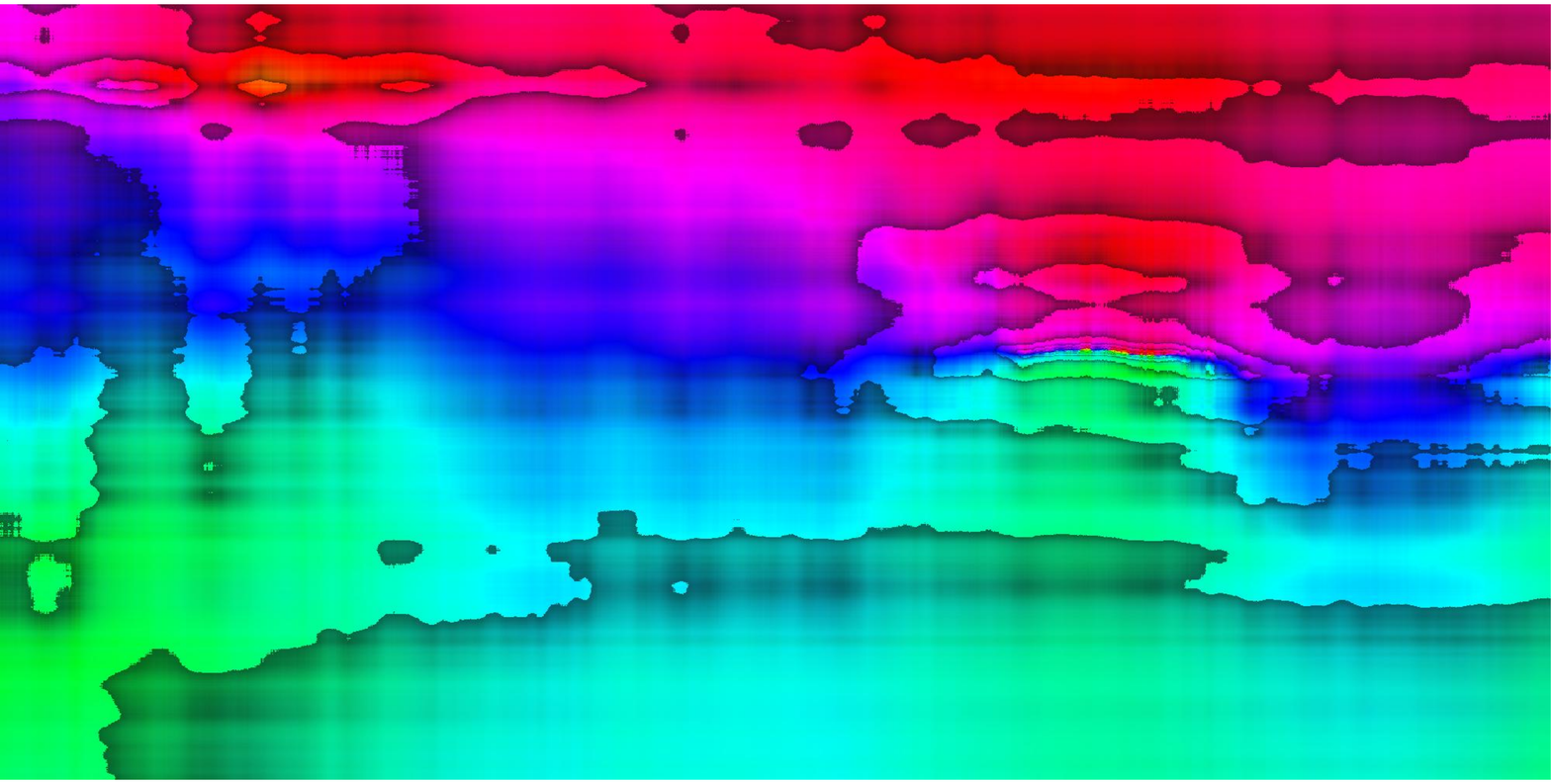}
		\caption{ A realization of the $S_{a,b}^{\Dord , \alpha}$ process with $\alpha=2$ (Gaussian distribution),  $a=1$, $b=-1$ and $\Dord = 1.3 - 0.7 \jj$.} \label{fig:Gaussian}
	\end{figure}

\begin{figure}[tb]
		\centering
	\includegraphics[width=12cm]{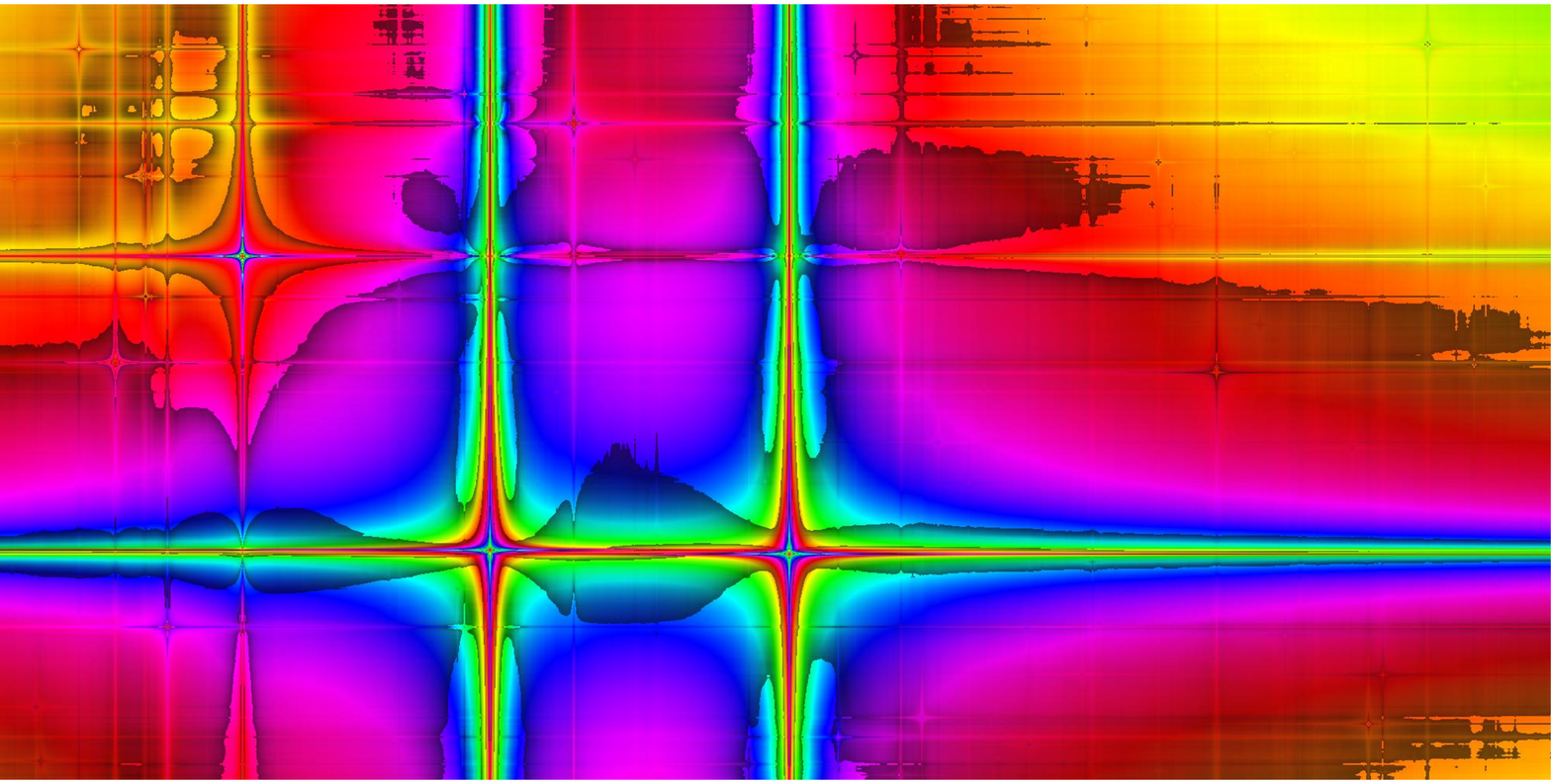}
		\caption{ A realization of the $S_{a,b}^{\Dord , \alpha}$ process with $\alpha=1$ (Cauchy distribution),  $a=1$, $b=1$ and $\Dord = 0.7 + \jj$.} \label{fig:Cauchy}
	\end{figure}

    The introduced $S_{a,b}^{\Dord , \alpha}$ processes in this work are complex-valued. Hence, to plot sample realizations we need to show the real and imaginary parts separately. Instead, for more compelling visual illustrations, we applied our framework to the generalizations of separable processes in $2$D 
     where complex values could be more conveniently shown by colors.
    In Figures \ref{fig:Gaussian} and \ref{fig:Cauchy} we can see two sample realizations. 
    For generating these figures, we have first generated a $2$D fine-grid discretization of the white noise process in form of an array of i.i.d. random variables. Then, we have applied the integration operator $\Iop_{a,b}^{\Dord ; k}$ to the white noise (the 2D array) both vertically and horizontally using the impulse response expressions in \eqref{eq:ImpulseResponseSimplified}. For Figure \ref{fig:Gaussian} we have used $\alpha=2$ (Gaussian distribution) with $a=1$, $b=-1$ and $\Dord = 1.3 - 0.7 \jj$; this case requires $k=\lfloor \frac{1}{2} + 1.3\rfloor = 1$. Similarly, we have set $\alpha=1$ (Cauchy distribution) with $a=b=1$ and $\Dord = 0.7 + \jj$ in Figure \ref{fig:Cauchy} (again requiring $k=\lfloor \frac{1}{1} + 0.7\rfloor = 1$). 
    The color-coding in these figures follows the standard approach for showing complex numbers, where the
    intensity of the pixels reflect the modulus of the complex numbers, while their phase are encoded in the hue.

\begin{figure}[tb]
		\centering
	\includegraphics[width=8cm]{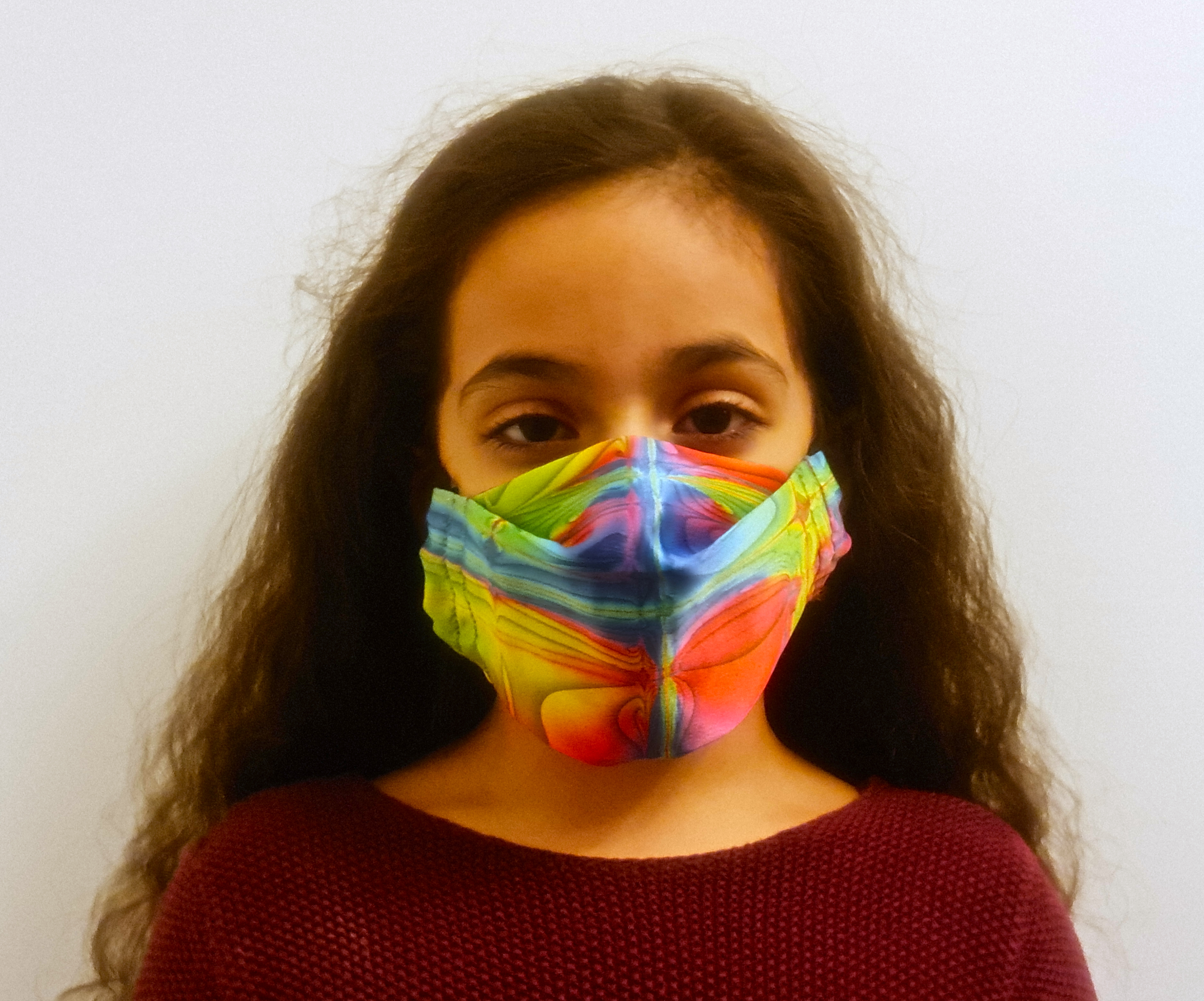}
		\caption{A face mask during the COVID-19 pandemic showing complex-order integration applied to a $2$D Poisson white noise.} \label{fig:Delaram}
	\end{figure}

It is interesting to mention that the $2$D plots of these processes were useful in designing face masks during the COVID-19 pandemic; in Figure \ref{fig:Delaram}, we see a face mask over which a fractionally integrated complex-valued process is printed. The white noise process in this case is Poisson with Gaussian jumps; although this process is not self-similar, it converges to the Gaussian white noise when the density of jumps goes to infinity~\cite{fageot2020gaussian}. Therefore, we expect the results generated by this Poisson white noise to resemble that of the Gaussian white noise with large enough jump densities.

\section{Useful Lemmas}\label{sec:Lemmas}

\begin{lemma}\label{lemma:L1property}
Let $h_{a,b}^{\Dord}(\cdot)$ be a homogeneous generalized function of degree $\Dord\in\C$ with $\real(\Dord)>-1$ as in \eqref{eq:hDef}, and let $n\in\{1,2,\dots,\lceil \real( \Dord)\rceil+1\}$. 
If $\widehat{\phi}:\R\mapsto\C$ is $n$-times continuously differentiable such that $|\widehat{\phi}^{(j)}(\omega)|(1+|\omega|^{r+j})$ is bounded for all $1\leq j\leq n$ and some $r >\real(\Dord) -n+1$, then, 
\begin{align*}
\forall\,1<T,  ~~~~~~~~ \Big\| \frac{\dd^n}{\dd \omega^n} \Big( \big(\widehat{\phi}(\omega)-\widehat{\phi}(\tfrac{\omega}{T})\big) h_{a,b}^{\Dord}(\omega)\Big) \Big\|_{1} 
\leq c \; T^{\real(\Dord) +2 -n},
\end{align*} 
where
\begin{align*}
c  = 
\tfrac{8 \, (r+1) \, (|\real(\Dord)|+2) \, \max(|a| \,,\, |b|)}{ (\real(\Dord) - n + 2) \, (r + n-\real(\Dord) - 1)} 
\bigg(\sum_{k=1}^{n} \scalebox{1}{$\binom{n}{k}$} \sup_{\omega} \Big(\big| \widehat{\phi}^{(k)}(\omega)\big| (1+|\omega|^{r+k})\Big) \bigg).
\end{align*}  
\end{lemma}

\textbf{Proof.}  
Intuitively, $\widehat{\phi}(\omega)-\widehat{\phi}(\tfrac{\omega}{T})$ behaves similar to $(1-\tfrac{1}{T})\omega\widehat{\phi}^{(1)}(\omega)$ around $\omega=0$ (a rigorous statement will be presented soon). Thus, $\big(\widehat{\phi}(\omega)-\widehat{\phi}(\tfrac{\omega}{T})\big) h_{a,b}^{\Dord}(\omega)$ has a zero of order at least $\real(\Dord) +1$ at $\omega=0$. This suggests that $\frac{\dd^n}{\dd \omega^n}  \big(\widehat{\phi}(\omega)-\widehat{\phi}(\tfrac{\omega}{T})\big)$ is either bounded at $\omega=0$ (if $n\leq \lceil \real(\Dord)\rceil$) or is at least locally integrable around $\omega=0$ (if $n= \lceil \real(\Dord)\rceil +1$). We also employ the boundedness property of $\big\{|\widehat{\phi}^{(j)}(\omega)|(1+|\omega|^{\gamma+j}) \big\}_i$ to show that $\frac{\dd^n}{\dd \omega^n}  \big(\widehat{\phi}(\omega)-\widehat{\phi}(\tfrac{\omega}{T})\big)$ decays asymptotically no slower than $\tfrac{1}{1+|\omega|^{1+\epsilon}}$ for some $\epsilon>0$ that is determined by $r$. As a result, we conclude that the $L_1$-norm of $\frac{\dd^n}{\dd \omega^n}  \big(\widehat{\phi}(\omega)-\widehat{\phi}(\tfrac{\omega}{T})\big)$ is finite ($\widehat{\phi}$ is not even required to be asymptotically decaying). However, the main point in Lemma \ref{lemma:L1property} is the scaling of the $L_1$-norm in terms of $T$.

To start our rigorous arguments, let us define
\begin{align}
\ushortw{\widehat{\phi}}(\omega)
= \sup_{\substack{\tau_{\phantom{.}} \\ |\tau|\geq |\omega|}} \big|\tfrac{\dd}{\dd \tau} \widehat{\phi}(\tau) \big|.
\end{align}
Based on the assumptions, $|\widehat{\phi}^{(1)}(\omega)|(1+|\omega|^{r+1})$ is bounded. Moreover, $r+1>\real(\Dord) - n + 2 > \real(\Dord) - \lceil \real(\Dord) \rceil +1 >0$. This shows that $|\widehat{\phi}^{(1)}(\omega)|$ is bounded and asymptotically decaying. Hence, $\ushortw{\widehat{\phi}}(\omega)$ is a bounded,  even, and non-negative-valued function that is non-increasing in terms of $|\omega|$. Since $r+1>0$, we conclude that $\ushortw{\widehat{\phi}}(\omega)(1+|\omega|^{r+1})$ is also bounded:
\begin{align}\label{eq:PhiBar}
\ushortw{\widehat{\phi}}(\omega)(1+|\omega|^{r+1})  
= \sup_{\substack{\tau_{\phantom{.}}  \\
|\tau|\geq |\omega|}} \big| \widehat{\phi}^{(1)}(\tau) \big| (1+|\omega|^{r+1}) 
&\leq \sup_{\substack{\tau_{\phantom{.}}  \\ 
|\tau|\geq |\omega|}} \big| \widehat{\phi}^{(1)}(\tau) \big| (1+|\tau|^{r+1}) \nonumber\\
&\leq \underbrace{ \sup_{\tau} \big| \widehat{\phi}^{(1)}(\tau) \big| (1+|\tau|^{r+1}) }_{\text{a finite constant}} .
\end{align}
Next, we bound the difference $\widehat{\phi}(\omega)-\widehat{\phi}(\tfrac{\omega}{T})$ via $\ushortw{\widehat{\phi}}(\omega)$. Because our approach is based on the Taylor's series, we need to initially separate the real and imaginary parts of $\widehat{\phi}$. Note that $\real\big(\widehat{\phi}(\omega)\big)$ and $\imag\big(\widehat{\phi}(\omega)\big)$ are both continuously differentiable functions for which we can write the Lagrange form of the Taylor's series as
\begin{align*}
\real\big(\widehat{\phi}(\omega)\big) - \real\big(\widehat{\phi}(\tfrac{\omega}{T})\big) &=  (\omega-\tfrac{\omega}{T}) \tfrac{\dd}{\dd \omega} \real \big(\widehat{\phi}(\omega)\big) \Big|_{\omega = \zeta_r} =   \tfrac{T-1}{T} \omega \,  \real \big(\widehat{\phi}^{(1)}(\zeta_r)\big) ,\\
\imag \big(\widehat{\phi}(\omega)\big) - \imag\big(\widehat{\phi}(\tfrac{\omega}{T})\big) &=  (\omega-\tfrac{\omega}{T}) \tfrac{\dd}{\dd \omega} \imag \big(\widehat{\phi}(\omega)\big) \Big|_{\omega = \zeta_i} =   \tfrac{T-1}{T} \omega \,  \imag \big(\widehat{\phi}^{(1)}(\zeta_i)\big) ,
\end{align*}
where $\zeta_r,\zeta_i$ are real numbers between $\tfrac{\omega}{T}$ and $\omega$. Therefore, if $\mathcal{I}_{\omega,T}$ represents the closed interval between $\frac{\omega}{T}$ and $\omega$, we have that
\begin{align}\label{eq:realPart}
\big| \real\big( \widehat{\phi}(\omega)-\widehat{\phi}(\tfrac{\omega}{T}) \big)\big| 
&\leq  \tfrac{T-1}{T} |\omega| \sup_{\tau \in \mathcal{I}_{\omega,T}} \big| \real \big( \phi^{(1)}(\tau) \big) \big|
\leq  |\omega| \sup_{\tau \in \mathcal{I}_{\omega,T}} \big|  \phi^{(1)}(\tau)\big| \nonumber\\
&\leq  |\omega|\; \ushortw{\widehat{\phi}}(\tfrac{\omega}{T})  .
\end{align}
Similarly, we can show that
\begin{align}\label{eq:imagPart}
\big| \imag\big( \widehat{\phi}(\omega)-\widehat{\phi}(\tfrac{\omega}{T}) \big)\big| 
\leq  |\omega|\;  \ushortw{\widehat{\phi}}(\tfrac{\omega}{T})  .
\end{align}
By combining \eqref{eq:realPart} and \eqref{eq:imagPart}, we achieve
\begin{align}
\big| \widehat{\phi}(\omega)-\widehat{\phi}(\tfrac{\omega}{T}) \big| 
\leq \big| \real\big( \widehat{\phi}(\omega)-\widehat{\phi}(\tfrac{\omega}{T}) \big)\big|  + \big| \imag\big( \widehat{\phi}(\omega)-\widehat{\phi}(\tfrac{\omega}{T}) \big)\big|  \leq 2|\omega|\;  \ushortw{\widehat{\phi}}(\tfrac{\omega}{T}) .
\end{align}
Now, we are equipped  to consider the main claim:
\begin{align*}
&\Big|\tfrac{\dd^n}{\dd \omega^n} \hspace{-0.5mm} \Big( \hspace{-1mm} \big(\widehat{\phi}(\omega)-\widehat{\phi}(\tfrac{\omega}{T})\big) h_{a,b}^{\Dord}(\omega) \hspace{-1mm}\Big) \hspace{-1mm} \Big| 
= \bigg| \hspace{-1mm} \sum_{k=0}^{n} \scalebox{1}{$\binom{n}{k}$} \Big(\tfrac{\dd^k}{\dd \omega^k} \big(\widehat{\phi}(\omega) \hspace{-1mm} -\hspace{-1mm} \widehat{\phi}(\tfrac{\omega}{T})\big) \hspace{-1mm} \Big) \hspace{-1mm} 
\Big(  \tfrac{\dd^{n-k}}{\dd \omega^{n-k}} h_{a,b}^{\Dord}(\omega) \hspace{-1mm}  \Big) \bigg|\nonumber\\
&~~~=  
\bigg|\big(\widehat{\phi}(\omega)-\widehat{\phi}(\tfrac{\omega}{T})\big) h_{a,b}^{\Dord,(n)}(\omega) + 
\sum_{k=1}^{n} \scalebox{1}{$\binom{n}{k}$} \big(\widehat{\phi}^{(k)}(\omega)-\tfrac{1}{T^k}\widehat{\phi}^{(k)}(\tfrac{\omega}{T})\big) h_{a,b}^{\Dord,(n-k)}(\omega)\bigg| \nonumber\\
&~~~ \leq  
2 \ushortw{\widehat{\phi}}(\tfrac{\omega}{T})  \; |\omega  \; h_{a,b}^{\Dord,(n)}(\omega)| +
\sum_{k=1}^{n} \scalebox{1}{$\binom{n}{k}$} \big(\big|\widehat{\phi}^{(k)}(\omega )\big| 
+ \tfrac{1}{T^{k}} \big|\widehat{\phi}^{(k)}(\tfrac{\omega}{T}) \big|\big) |h_{a,b}^{\Dord,(n-k)}(\omega)|
\end{align*}
which implies
\begin{align}\label{eq:L1bound1}
& \Big\|\tfrac{\dd^n}{\dd \omega^n} \Big( \big(\widehat{\phi}(\omega)-\widehat{\phi}(\tfrac{\omega}{T})\big) h_{a,b}^{\Dord}(\omega)\Big) \Big\|_{1} 
 \leq 2\big\|\ushortw{\widehat{\phi}}(\tfrac{\omega}{T}) \; \omega  \; h_{a,b}^{\Dord,(n)}(\omega)\big\|_{1}
 \nonumber\\
& \hspace{5em} +\sum_{k=1}^{n} \scalebox{1}{$\binom{n}{k}$} \big(\big\|\widehat{\phi}^{(k)}(\omega) h_{a,b}^{\Dord,(n-k)}(\omega)\big\|_{1} 
+ \tfrac{1}{T^{k}} \big\|\widehat{\phi}^{(k)}(\tfrac{\omega}{T}) h_{a,b}^{\Dord,(n-k)}(\omega)\big\|_{1}\big) . 
\end{align}
To simplify the upperbound, note that if $l(\omega)$ is a homogeneous generalized function of degree $s$, then, for a generic function $g(\omega)$ we have that
\begin{align*}
\big\| g(\tfrac{\omega}{T}) l(\omega) \big\|_1 
&= \int_{\R} \big| g(\underbrace{\tfrac{\omega}{T}}_{\nu} ) l(\omega)\big| \dd \omega = T \int_{\R} \big| g(\nu ) l(T\,\nu)\big| \dd \nu \nonumber\\
&= T^{\real(s) + 1} \int_{\R} \big| g(\nu ) l(\nu)\big| \dd \nu = T^{\real(s) + 1} \big\| g(\omega) l(\omega) \big\|_1.
\end{align*}
Since $\omega \,h_{a,b}^{\Dord,(n)}(\omega)$ and $h_{a,b}^{\Dord,(n-k)}(\omega)$ are both homogeneous with degrees $\Dord-n+1$ and $\Dord-n+k$, respectively, we can simplify \eqref{eq:L1bound1} as
\begin{align}
\Big\|\tfrac{\dd^n}{\dd \omega^n} \Big( \big(&\widehat{\phi}(\omega) -\widehat{\phi}(\tfrac{\omega}{T})\big) h_{a,b}^{\Dord}(\omega)\Big) \Big\|_{1} 
 \leq \phantom{+}
 2T^{\real(\Dord)+2-n} \big\|\omega \; \ushortw{\widehat{\phi}}(\omega) \; h_{a,b}^{\Dord,(n)}(\omega)\big\|_{1} 
\nonumber\\
&+ \big(T^{\real(\Dord)+1-n}+1\big)\sum_{k=1}^{n} \scalebox{1}{$\binom{n}{k}$} \Big\|\widehat{\phi}^{(k)}(\omega) h_{a,b}^{\Dord,(n-k)}(\omega)\Big\|_{1} \nonumber\\
< & \underbrace{ \Big(2\big\|\omega \; \ushortw{\widehat{\phi}}(\omega) \; h_{a,b}^{\Dord,(n)}(\omega)\big\|_{1} +   2\sum_{k=1}^{n} \scalebox{1}{$\binom{n}{k}$} \Big\|\widehat{\phi}^{(k)}(\omega) h_{a,b}^{\Dord,(n-k)}(\omega)\Big\|_{1} \Big)}_{ c } T^{\real(\Dord)+2-n}
\end{align}
which proves the claim. However, we still need to show that $\big\|\omega \; \ushortw{\widehat{\phi}}(\omega) \; h_{a,b}^{\Dord,(n)}(\omega)\big\|_{1} $ and $\big\|\widehat{\phi}^{(k)}(\omega) h_{a,b}^{\Dord,(n-k)}(\omega)\big\|_{1}$ are all well-defined and finite. To that end, we first observe that 
\begin{align*}
\forall\, 0\leq k\leq n:~~\big|h_{a,b}^{\Dord,(k)}(\omega)\big| \leq \ab \, |\real(\Dord) \hspace{-1mm}-\hspace{-1mm} k|\, \cdot \, |\omega|^{\real(\Dord)-k},
\end{align*}
where $\ab=\max\big(|a|\,,\,|b|\big)$. 
Based on the assumption, the value $\widehat{\phi}_k = \sup_{\omega} \,\big|\widehat{\phi}^{(k)}(\omega) \big| (1+|\omega|^{r+k})$ is finite for all $1\leq k\leq n$. Moreover, we have previously shown in \eqref{eq:PhiBar} that $\ushortw{\widehat{\phi}}(\omega)(1+|\omega|^{r+1})$ is upper-bounded by $\widehat{\phi}_1$. 
Therefore,
\begin{align*}
\big\|\omega \; \ushortw{\widehat{\phi}}(\omega) \; &h_{a,b}^{\Dord,(n)}(\omega)\big\|_{1} 
  \leq 2 |\real(\Dord) \hspace{-1mm}-\hspace{-1mm} n|\, \ab \, \widehat{\phi}_1 \int_{0}^{\infty}  \frac{\omega^{\real(\Dord)-n+1}}{ 1+\omega^{r+1} }\dd\omega \nonumber\\
 & \leq 
 2 |\real(\Dord) \hspace{-1mm}-\hspace{-1mm} n|\, \ab \, \widehat{\phi}_1 \bigg( \int_{0}^{1}  \omega^{\real(\Dord)-n+1} \dd\omega 
 +   \int_{1}^{\infty}  \omega^{\real(\Dord)-n-r}\dd\omega \bigg)\nonumber\\
  & = 
 2 |\real(\Dord) \hspace{-1mm}-\hspace{-1mm} n|\, \ab \, \widehat{\phi}_1 \bigg( \frac{ \omega^{\real(\Dord)-n+2} }{ \real(\Dord)-n+2 } \bigg|_0^1 
 +   \frac{ \omega^{\real(\Dord)-n-r +1} }{ \real(\Dord)-n-r +1 } \bigg|_{1}^{\infty} \bigg).
\end{align*}
The upper-bound is finite because
\begin{eqnarray}\label{eq:finiteness}
& \real(\Dord)-n+2 \geq \real(\Dord)-(\lceil \real(\Dord)\rceil + 1) +2 = \real(\Dord)-\lceil \real(\Dord)\rceil + 1 >0, &\nonumber\\
& \real(\Dord)-n-r +1 < \real(\Dord)-n-(\real(\Dord) - n + 1) +1 = 0. &
\end{eqnarray}
For the existence of $\big\|\widehat{\phi}^{(k)}(\omega) h_{a,b}^{\Dord,(n-k)}(\omega)\big\|_{1}$, we employ a similar technique:
\begin{align*}
\big\|\widehat{\phi}^{(k)}&(\omega)  h_{a,b}^{\Dord,(n-k)}(\omega)\big\|_{1}
  \leq 2 |\real(\Dord) \hspace{-1mm}-\hspace{-1mm} n \hspace{-1mm}+\hspace{-1mm} k|\, \ab \, \widehat{\phi}_k \int_{0}^{\infty}  \frac{\omega^{\real(\Dord)-n+k}}{ 1+\omega^{r+k} }\dd\omega \nonumber\\
 & \leq 
 2|\real(\Dord) \hspace{-1mm}-\hspace{-1mm} n \hspace{-1mm}+\hspace{-1mm} k|\, \ab \, \widehat{\phi}_k \bigg( \int_{0}^{1}  \omega^{\real(\Dord)-n+k} \dd\omega 
 +   \int_{1}^{\infty}  \omega^{\real(\Dord)-n-r}\dd\omega \bigg)\nonumber\\
  & = 
 2 |\real(\Dord) \hspace{-1mm}-\hspace{-1mm} n \hspace{-1mm}+\hspace{-1mm} k|\, \ab\, \widehat{\phi}_k \bigg( \frac{ \omega^{\real(\Dord)-n+k+1} }{ \real(\Dord)-n+k+1 } \bigg|_0^1 
 +   \frac{ \omega^{\real(\Dord)-n-r +1} }{ \real(\Dord)-n-r +1 } \bigg|_{1}^{\infty} \bigg).
\end{align*}
The finiteness of the upper-bound follows from \eqref{eq:finiteness} by considering $k\geq 1$.  $\hspace{\stretch{1}}\blacksquare$



\begin{lemma}\label{lemma:ArashTechnique}
If  $\phi(x):\R\setminus\{0\}\mapsto\C$ satisfies
\begin{align*}
\forall\, x\in\R\setminus\{0\},~ \forall\, T\in[2,4[:~~~ \big|\phi(x) - T^{\Dord}\phi(T\,x)\big| \leq \frac{c}{|x|^m},
\end{align*}
where $c\in\R^{+}$, $m\in\R$ and $\Dord\in\C$ are constants, then, there is a constant $\bar{c}$ such that
\begin{align*}
\forall\, x\in\R\setminus\{0\} :~~~ \big| \phi(x)\big| \leq \frac{\bar{c} }{\min\big(|x|^{\real(\Dord)} \,,\, |x|^{m} \big)},
\end{align*}
if $m\neq\real(\Dord)$, and 
\begin{align*}
\forall\, x\in\R\setminus\{0\} :~~~ \big| \phi(x)\big| \leq \bar{c}\;\frac{1+\big|\log|x|\big|}{|x|^{\real(\Dord)}},
\end{align*}
if $m=\real(\Dord)$.

\end{lemma}

\textbf{Proof.}  
Let $\tau$ be an arbitrary positive real. Any real number $x$ with $2\tau\leq|x|$ or $0<|x| \leq \tfrac{\tau}{2}$ can be uniquely written in the form of $\pm T^{n}\,\tau$ with $n\in\{\pm 1,\pm2,\dots,\pm2^i,\dots\}$ and $T\in[2,4[$\footnote{This can be achieved by setting $n={\rm sign}(\theta) \, 2^{\lfloor \log_2|\theta|\rfloor}$ and $T= 2^{2^{\{\log_2 \theta\}}}$, where $\theta =\log_2\big|\tfrac{x}{\tau}\big|$.}. For the sake of simplicity, we continue the proof only for $0<x=\tau T^{n}$ (the negative case is similar). We define
$$S_n = \left\{ \begin{array}{ll}
\{0,1,\dots,n-1\}, & n\geq 1, \phantom{\big|} \\
\{-1,-2,\dots,n\}, & n\leq -1. \phantom{\big|}
\end{array}\right.$$
Now, we have that
\begin{align}\label{eq:BeforeTwoCases}
\Big| \tau^{\Dord} \phi(\tau) - \big(T^{n}\tau\big)^{\Dord} \phi(T^{n}\,\tau) \Big| 
&\leq \sum_{k\in S_n} \Big| \big(T^{k}\tau\big)^{\Dord} \phi(T^{k}\,\tau) - \big(T^{k+1}\tau\big)^{\Dord} \phi( T^{k+1}\,\tau) \Big| \nonumber\\
& = \sum_{k\in S_n} \big(T^{k}\tau\big)^{\real(\Dord)} \Big|  \phi(T^{k}\,\tau) - T^{\Dord} \phi( T^{k+1}\,\tau) \Big| \nonumber\\
& \leq \sum_{k\in S_n} \frac{ c \big(T^{k}\tau\big)^{\real(\Dord)}}{ (T^{k}\,\tau)^m } 
=  c \; \tau^{\real(\Dord) - m}\sum_{k\in S_n} \big(T^{\real(\Dord)-m}\big)^k .
\end{align}
We consider the two cases of $m\neq\real(\Dord)$ and $m=\real(\Dord)$ separately.

\begin{enumerate}
\item $m\neq \real(\Dord)$. Thus, $\sum_{k\in S_n} T^{k(\real(\Dord) - m)} = \big| \tfrac{T^{n(\real(\Dord)-m)} - 1}{T^{\real(\Dord)-m} - 1} \big|$ and \eqref{eq:BeforeTwoCases} can be written as
\begin{align}\label{eq:mainInEq}
&\Big| \tau^{\Dord} \phi(\tau) - \big(T^{n}\tau\big)^{\Dord} \phi(T^{n}\,\tau) \Big| 
\leq \tfrac{c  }{ |T^{\real(\Dord)-m} - 1| } \big(T^n\tau\big)^{\real(\Dord)-m} + \tfrac{ c \; \tau^{\real(\Dord) - m} }{ |T^{\real(\Dord)-m} - 1 |} \nonumber\\
& \Rightarrow ~~~ \Big|\big(T^{n}\tau\big)^{\Dord} \phi(T^{n}\,\tau) \Big| \leq \tfrac{c  (T^n\tau)^{\real(\Dord)-m}  }{ |T^{\real(\Dord)-m} - 1| } + \tfrac{ c \; \tau^{\real(\Dord) - m} }{ |T^{\real(\Dord)-m} - 1 |} +  \big|\tau^{\Dord} \phi(\tau)\big|  \nonumber\\
 & \leq \big( c\tfrac{ \tau^{\real(\Dord) - m}+1 }{ |T^{\real(\Dord)-m} - 1 |} +  \big|\tau^{\Dord} \phi(\tau)\big|\big) \; \max\big((T^n\tau)^{\real(\Dord)-m} \,,\, 1\big) \nonumber\\
  & \leq \underbrace{\big( c\tfrac{ \tau^{\real(\Dord) - m}+1 }{ \min(|2^{\real(\Dord)-m} - 1 | \,,\, |4^{\real(\Dord)-m} - 1 |)} +  \big|\tau^{\Dord} \phi(\tau)\big|\big)}_{\bar{c}_{\tau}} \; \max\big((T^n\tau)^{\real(\Dord)-m} \,,\, 1\big),
\end{align}
where $\bar{c}_{\tau}$ is a constant that depends  neither on $n$ nor on $T$. As explained earlier, for any positive $x\not\in \,]\tfrac{\tau}{2},2\tau[$ we can find suitable $n$ and $T$, such that $x= T^n\tau$. By substituting $x= T^n\tau$ in \eqref{eq:mainInEq} we obtain
\begin{align*}
\forall \, 0<x\not\in \, ]\tfrac{\tau}{2} , 2\tau[:~~~ \big|\phi(x) \big| \leq \frac{\bar{c}_{\tau}}{ \min\big(x^{\real(\Dord)} \,,\, x^{m}\big) }.
\end{align*}

\item $m=\real \, z$. This implies that $\sum_{k\in S_n} T^{k(\real(\Dord) - m)} = |n|$. Now \eqref{eq:BeforeTwoCases} simplifies to
\begin{align}
 &\Big| \tau^{\Dord} \phi(\tau) - \big(T^{n}\tau\big)^{\Dord} \phi(T^{n}\,\tau) \Big| 
\leq c \; \tau^{\real(\Dord) - m}\; |n| = c \; \tau^{\real(\Dord) - m}\; \Big| \tfrac{\log(T^{n}\tau)-\log(\tau)}{\log T}\Big| \nonumber\\
&\Rightarrow ~  \Big| \big(T^{n}\tau\big)^{\Dord} \phi(T^{n}\,\tau) \Big| 
 \leq \underbrace{ \big( c\tfrac{ \tau^{\real(\Dord) - m}(1+|\log \tau|)}{\log 2} + | \tau^{\Dord} \phi(\tau)|\big)}_{\bar{c}_{\tau}} \big( 1 \hspace{-0.5mm} + \hspace{-0.5mm} |\log(T^{n}\,\tau)| \big),
\end{align}
where $\bar{c}_{\tau}$ is again a constant that depends neither on $n$ nor on $T$. Thus, similar to the previous case, we conclude that
\begin{align*}
\forall \, 0<x\not\in \, ]\tfrac{\tau}{2} , 2\tau[:~~~ \big|\phi(x) \big| \leq \bar{c}_{\tau}\,\frac{1+|\log x|}{ x^{\real(\Dord)}  }.
\end{align*}
\end{enumerate}

The above arguments prove the claim, except for $x\in \, ]\tfrac{\tau}{2} , 2\tau[$. However, the choice of $\tau$ was arbitrary; thus, all the above results also hold for $\tau'=4\tau$, if $\bar{c}_{\tau}$ is replaced with $\bar{c}_{\tau'}$. In addition, $]\tfrac{\tau}{2} , 2\tau[$ and $]\tfrac{\tau'}{2} , 2\tau'[$ have empty intersection. Thus, if we use $\bar{c}=\max(\bar{c}_{\tau} , \bar{c}_{\tau'})$ as the constant, the inequalities hold for all $x\neq 0$. $\hspace{\stretch{1}}\blacksquare$

\begin{lemma}\label{lemma:ContDeriv}
If $\varphi:\R \rightarrow \C$ is $k+n$ times continuously differentiable ($k\geq 1$ and $n\geq 0$), then, the function $\phi(x)=\frac{\varphi(x) - \sum_{j=0}^{k-1} \tfrac{x^j}{j!}\varphi^{(j)}(0)}{x^{k}}$ is $n$ times continuously differentiable. In addition, if $\tfrac{|\varphi^{(j)}(x)|}{1+|x|^{k-1-j}}$ is bounded for all $0\leq j< k$, then, $|\phi^{(n)}(x)| (1+|x|^{n+1})$ is also bounded.
\end{lemma}

\textbf{Proof.}  
For notational convenience, we define
\begin{align*}
f(x) = \sum_{i=0}^{k-1} \tfrac{x^i}{i!}\varphi^{(i)}(0) ~~~~~~,~~~~~~
g(x) = \sum_{i=0}^{k+n-1} \tfrac{x^i}{i!}\varphi^{(i)}(0).
\end{align*}
Since $\varphi$ is $k+n$ times continuously differentiable, $\varphi^{(j)}$ is also $k+n-j$ times continuously differentiable. 
By applying the Lagrange's form of the Taylor series on the real and imaginary parts of $\varphi^{(j)}(x)$ separately for $0\leq j\leq n$, and then combining them, we have that
\begin{align}\label{eq:TaylorForms}
\varphi^{(j)}(x) &= \hspace{-2mm} \sum_{i=0}^{k+n-j-1} \hspace{-2mm} \tfrac{x^i}{i!}\varphi^{(j+i)}(0) \hspace{-1mm}+\hspace{-1mm} \tfrac{x^{k+n-j}}{(k+n-j)!} \hspace{-0.5mm} \Big( \hspace{-1mm}\real \hspace{-1mm} \big(\varphi^{(k+n)}(\zeta_{x,R}^{(j)})\big)  \hspace{-1.5mm}+\hspace{-1.0mm} \jj \imag \hspace{-1mm} \big(\varphi^{(k+n)}(\zeta_{x,I}^{(j)})\big) \hspace{-1mm} \Big) \nonumber\\
&= \tfrac{\dd^{j}}{\dd x^j} g(x) \hspace{-1mm}+\hspace{-1mm} \tfrac{x^{k+n-j}}{(k+n-j)!} \hspace{-0.5mm} \Big( \hspace{-1mm}\real \hspace{-1mm} \big(\varphi^{(k+n)}(\zeta_{x,R}^{(j)})\big)  \hspace{-1.5mm}+\hspace{-1.0mm} \jj \imag \hspace{-1mm} \big(\varphi^{(k+n)}(\zeta_{x,I}^{(j)})\big) \hspace{-1mm} \Big),
\end{align}
where $\zeta_{x,R}^{(j)}$ and $\zeta_{x,I}^{(j)}$ are real numbers between $0$ and $x$. Let us now evaluate the $m$th ($0\leq m\leq n$) derivative of ${\phi}$:
\begin{align}\label{eq:SimpleDerivForm}
\tfrac{\dd^m}{\dd x^m} {\phi}(x) 
&= \tfrac{\dd^m}{\dd x^m} \Big(\frac{{\varphi}(x) -{f}(x)}{x^{k}} \Big)
=\sum_{j=0}^{m} \scalebox{1}{$\binom{m}{j}$} \big( {\varphi}^{(j)}(x) -{f}^{(j)}(x) \big) ~ \tfrac{\dd^{m-j}}{\dd x^{m-j}}\big(\tfrac{1}{x^{k}}\big).
\end{align}
As $\{ {\varphi}^{(j)} \}_j$ are continuous for all $0\leq j\leq m$, it is straightforward to see that ${\phi}^{(m)}(x)$ exists and is continuous at $x \neq 0$. To investigate the behaviour of ${\phi}^{(m)}(x)$ around $x=0$, we rewrite \eqref{eq:SimpleDerivForm} by using \eqref{eq:TaylorForms} as
\begin{align*}
&\tfrac{\dd^m}{\dd x^m} {\phi}(x) 
= \sum_{j=0}^{m} \scalebox{1}{$\binom{m}{j}$} \bigg(  {g}^{(j)}(x) -{f}^{(j)}(x) + \tfrac{x^{k+n-j}}{(k+n-j)!} 
\Big(\real\big( {\varphi}^{(k+n)}(\zeta_{x,R}^{(j)}) \big)  \nonumber\\
& \phantom{\tfrac{\dd^m}{\dd x^m} {\phi}(x) 
=\sum_{j=0}^{m} \scalebox{1}{$\binom{m}{j}$} \bigg(  {g}^{(j)}(x) -{f}^{(j)}(x) +}
+ \jj \imag\big( {\varphi}^{(k+n)}(\zeta_{x,I}^{(j)}) \big) \Big)
  \bigg)  \tfrac{\dd^{m-j}}{\dd x^{m-j}}\Big(\tfrac{1}{x^{k}}\Big) \nonumber\\
&= \sum_{j=0}^{m} \scalebox{1}{$\binom{m}{j}$} \big(  {g}^{(j)}(x) -{f}^{(j)}(x) \big) ~ \tfrac{\dd^{m-j}}{\dd x^{m-j}}\big(\tfrac{1}{x^{k}}\big) 
\nonumber\\
& ~~ +
 \sum_{j=0}^{m} \scalebox{1}{$\binom{m}{j}$} \tfrac{x^{k+n-j}}{(k+n-j)!}
\Big( \hspace{-1.5mm} \real \hspace{-0.8mm} \big( {\varphi}^{(k+n)}(\zeta_{x,R}^{(j)}) \hspace{-0.5mm} \big) \hspace{-1mm}+\hspace{-1mm} \jj \imag \hspace{-0.8mm}  \big( {\varphi}^{(k+n)}(\zeta_{x,I}^{(j)}) \hspace{-0.5mm} \big) \hspace{-1.5mm}  \Big) 
 \tfrac{\dd^{m-j}}{\dd x^{m-j}}\big(\tfrac{1}{x^{k}}\big) \nonumber\\
&= 
\tfrac{\dd^m}{\dd x^m} \big(\underbrace{\tfrac{{g}(x) - {f}(x)}{w^{k}}}_{{P}(x)}\big) + \tfrac{x^{n-m}}{(k-1)! \, (n+ 1 - m)!} \sum_{j=0}^{m}  \tfrac{(-1)^{m-j} \binom{m}{j} }{\binom{k+n-j}{k+m-j-1}} \; 
\Big( \hspace{-1mm} \real\big( {\varphi}^{(k+n)}(\zeta_{x,R}^{(j)}) \big) \nonumber\\
&\phantom{\tfrac{\dd^m}{\dd x^m} \big(\underbrace{\tfrac{{g}(x) - {f}(x)}{x^{k}}}_{{P}(x)}\big) + \tfrac{x^{n-m}}{(k-1)! \, (n+ 1 - m)!} \sum_{j=0}^{m}  \tfrac{(-1)^{m-j} \binom{m}{j} }{\binom{k+n-j}{k+m-j-1}} \; 
\Big(}
~~ + \jj \imag\big( {\varphi}^{(k+n)}(\zeta_{x,I}^{(j)}) \big) \hspace{-1mm} \Big) ,
\end{align*}
where ${P}(x)$ is a polynomial of degree no more than $n-1$. 
Since $\{\zeta_{x,R}^{(j)} , \zeta_{x,I}^{(j)}\}_j$ are all between $0$ and $x$, they all approach $0$ when $x \rightarrow 0$. Furthermore, ${\varphi}^{(k+n)}$ is a continuous function by assumption. Thus
\begin{align*}
\lim_{x \rightarrow 0} \tfrac{\dd^m}{\dd x^m} {\phi}(x) 
&=  \left\{ \begin{array}{ll}
{P}^{(m)}(0), & 0\leq m <n,  \phantom{\Big|}\\
\tfrac{n!}{(k+n)!} \; {\varphi}^{(k+n)}(0), & m=n. \phantom{\Big|}
\end{array}\right.
\end{align*}
Now that we have shown ${\phi}^{(n)}(x) $ is continuous, we turn to the claimed decay result. Because of the continuity of ${\phi}^{(n)}$, it is sufficient to show that $|{\phi}^{(n)}(x)| (1+|x|^{n+1})$ is bounded only for $|x|\geq 1$. 
As ${f}^{(j)}(x)$ is a polynomial of degree at most $\max(0 \,,\, k-j-1)$, it is evident that $\frac{ {f}^{(j)}(x) }{1+|x|^{k-j-1}}$ is bounded. This result, in addition to the assumption in the lemma, implies that 
$$\frac{{\varphi}^{(j)}(x) - {f}^{(j)}(x) }{1+|x|^{k-j-1}}$$
is also bounded. Further, we have that
\begin{align*}
\frac{(1+|x|^{n+1})(1+|x|^{k-j-1})}{|x|^{k+n-j}} 
\stackrel{|x|\geq 1}{\leq} 
\frac{(|x|^{n+1}+|x|^{n+1})(|x|^{k-j-1}+|x|^{k-j-1})}{|x|^{k+n-j}} = 4.
\end{align*}
Now, recalling \eqref{eq:SimpleDerivForm}, we  obtain
\begin{align*}
 |{\phi}^{(n)}(x) | (1+|x|^{n+1})
\leq &
 (1+|x|^{n+1}) \sum_{j=0}^{n} \scalebox{1}{$\binom{n}{j}$} \tfrac{(k+n-j-1)!}{(k-1)! \, |x|^{k+n-j}} \, \big| {\varphi}^{(j)}(x) -{f}^{(j)}(x) \big| 
 \nonumber\\
 = & \sum_{j=0}^{n} \tfrac{\binom{n}{j} \, (k+n-j-1)!}{(k-1)! } \tfrac{| {\varphi}^{(j)}(x) -{f}^{(j)}(x) |}{ 1+|x|^{k-j-1} } \tfrac{(1+|x|^{n+1})(1+|x|^{k-j-1})}{|x|^{k+n-j}}.
\end{align*}
The upper-bound consists of finitely many bounded terms for $|x|\geq 1$, and is therefore, bounded. 
$\hspace{\stretch{1}} \blacksquare$


\begin{corr}\label{corr:SimpleIntegration}
For $\varphi\in\mathcal{S}$, we know that $\widehat{\varphi}$ (Fourier transform of $\varphi$) is infinitely differentiable and $\tfrac{|\widehat{\varphi}^{(j)}(\omega)|}{1+|\omega|^{m}}$ is bounded for all $0\leq j,m$. Using  Lemma \ref{lemma:ContDeriv}, we conclude that $\widehat{\phi}(\omega)=\frac{\widehat{\varphi}(\omega) - \sum_{i=0}^k \tfrac{\omega^i}{i!}\widehat{\varphi}^{(i)}(0)}{\omega^{k+1}}$ is infinitely differentiable and $\widehat{\phi}^{(j)}\in L_1$ for all $j\geq 1$. This implies that $(1+|x|^m)\phi(x)$ (where $\phi$ is the inverse Fourier of $\widehat{\phi}$) is bounded for all $m\in\N$. Nevertheless, since $\widehat{\phi}\not\in L_1$, $\phi(x)$ is not continuous.
\end{corr}


\begin{lemma}\label{lemma:BasicsIntegral}
Let $h_{a,b}^{\Dord}(\cdot)$ be a homogeneous function of degree $\Dord\in\C$  as in \eqref{eq:hDef} and let $k\geq \max\big(1 \,,\, \lfloor \real (\Dord)\rfloor  \big)$. 
For a $k$-times continuously differentiable $\widehat{\varphi}$, define
\begin{align*}
\widehat{\varphi}_{a,b}^{\Dord; \scalebox{0.6}{$k$}}(\omega) = \frac{\widehat{\varphi}(\omega)-\sum_{i=0}^{k-1}\tfrac{\widehat{\varphi}^{(i)}(0)}{i!}\omega^i }{h_{a,b}^{\Dord}(\omega)}.
\end{align*}
If $\tfrac{|\widehat{\varphi}(\omega)|}{1+|\omega|^{k-1}}$ and $\Big\{\tfrac{|\widehat{\varphi}^{(j)}(\omega)|}{1+|\omega|^{k-j}} \Big\}_{j=1}^{k}$ are all bounded, then, the inverse Fourier of $\widehat{\varphi}_{a,b}^{\Dord; \scalebox{0.6}{$k$}}$ denoted by ${\varphi}_{a,b}^{\Dord; \scalebox{0.6}{$k$}}(x)$ exists as a function for $x\neq 0$, and 
\begin{enumerate}[(i)]
\item \label{claim:k_1_NonZ} when $\real (\Dord)\not\in \Z$  and $k= \lfloor \real (\Dord)\rfloor $, the function ${\varphi}_{a,b}^{\Dord; \scalebox{0.6}{$k$}}(x)$ is continuous everywhere including $x=0$.

\item \label{claim:k_1_Z} when $\real (\Dord)\in \Z$, $\imag (\Dord)\neq 0$ and $k= \lfloor \real (\Dord)\rfloor $, the function ${\varphi}_{a,b}^{\Dord; \scalebox{0.6}{$k$}}(x)$ is continuous at $x\neq 0$ and bounded around $x=0$.

\item \label{claim:k_1_ZZh} when $\Dord\in \Z$, $h_{a,b}^{\Dord}(\omega)\equiv c\, \omega ^{\Dord}$ for some $c\in\C$, and $k= \lfloor \real (\Dord)\rfloor $, the function ${\varphi}_{a,b}^{\Dord; \scalebox{0.6}{$k$}}(x)$ is continuous at $x\neq 0$ and bounded around $x=0$.

\item \label{claim:k_1_ZZnh} when $\Dord\in \Z$, $h_{a,b}^{\Dord}(\omega)\not\equiv c\,  \omega ^{\Dord}$, and $k= \lfloor \real (\Dord)\rfloor $, the function ${\varphi}_{a,b}^{\Dord; \scalebox{0.6}{$k$}}(x)$ is continuous at $x\neq 0$, and $\tfrac{|\varphi_{a,b}^{\Dord; \scalebox{0.6}{$k$}}(x)|}{\log|x|}$ is bounded around $x=0$.

\item \label{claim:k_NonZ} when $k> \lfloor \real (\Dord)\rfloor $, the function ${\varphi}_{a,b}^{\Dord; \scalebox{0.6}{$k$}}(x)$ is continuous at $x\neq 0$, and $|x|^{k-\real(\Dord) } |\varphi_{a,b}^{\Dord; \scalebox{0.6}{$k$}}(x)|$ is bounded at $x\in [-1,1]$.
\end{enumerate}
\end{lemma}

\textbf{Proof.}
We first rewrite $\widehat{\varphi}_{a,b}^{\Dord; \scalebox{0.6}{$k$}}$ as
\begin{align}
\widehat{\varphi}_{a,b}^{\Dord; \scalebox{0.6}{$k$}}(\omega) = \underbrace{\frac{\widehat{\varphi}(\omega)-\sum_{i=0}^{k-1}\tfrac{\widehat{\varphi}^{(i)}(0)}{i!}\omega^i}{\omega^{k}} }_{\widehat{\phi}(\omega)} \underbrace{ \frac{\omega^{k}}{h_{a,b}^{\Dord}(\omega)} }_{\ut{\scalebox{0.8}{$h$}}(\omega)} = \widehat{\phi}(\omega) {\ut{$h$}}(\omega).
\end{align}
Lemma \ref{lemma:ContDeriv} implies that $\widehat{\phi}$ is continuous and $\widehat{\phi}(\omega)(1+|\omega|)$ is bounded. In addition,  ${\ut{$h$}}(\omega) = h_{a',b'}^{\Dord'}(\omega)$ with $a'=\frac{1}{a},~b'=\frac{(-1)^k}{b},~\Dord'=k-\Dord$ is a homogeneous function of degree $\ut{$\Dord$}=k-\Dord$. 
To simplify, let $\lambda=\real (\ut{$\Dord$})$; one can check that $-1<\lambda$. As $|{\ut{$h$}}(\omega)|$ is proportional to $|\omega|^{\lambda}$, we know that ${\ut{$h$}}(\omega)$ is locally integrable ($-1<\lambda$). Here, we continue as
\begin{align}\label{eq:PhiH}
\widehat{\varphi}_{a,b}^{\Dord; \scalebox{0.6}{$k$}}(\omega) =& \Big( \underbrace{ \widehat{\phi}(\omega) + \tfrac{ \widehat{\varphi}^{(k-1)}(0) }{ (k-1) ! \, \omega } \indc_{|\omega|>1} }_{ \widehat{ \ut{$\phi$} }(\omega) }
- \tfrac{ \widehat{\varphi}^{(k-1)}(0) }{ (k-1) ! \, \omega } \indc_{|\omega|>1} \Big) {\ut{$h$}}(\omega) \nonumber\\
=& \widehat{ \ut{$\phi$} }(\omega) \widehat{\ut{$h$}}(\omega) 
- \tfrac{ \widehat{\varphi}^{(k-1)}(0) }{ (k-1) ! } |\omega|^{\lambda-1+\jj\imag(\Dord)} \big(\tfrac{\indc_{\omega>1}}{a} + \tfrac{(-1)^{k-1}\indc_{\omega<-1} }{b} \big),
\end{align}
where $\indc_{I}$ is the indicator function for the set $I$. The function $\widehat{ \ut{$\phi$} }(\omega)$ is essentially equal to $\frac{\widehat{\varphi}(\omega)-\sum_{i=0}^{k-2}\tfrac{\widehat{\varphi}^{(i)}(0)}{i!}\omega^i}{\omega^{k}}$ for $|\omega|>1$, and is bounded for $|\omega|<1$. By invoking Lemma \ref{lemma:ContDeriv} it is not difficult to see that $|\widehat{ \ut{$\phi$} }(\omega)|(1+\omega^2)$ is bounded. Thus, $|\widehat{ \ut{$\phi$} }(\omega)\widehat{\ut{$h$}}(\omega)| \leq c\tfrac{|\omega|^{\lambda}}{1+\omega^{2}}\in L_1$; this implies that $\widehat{ \ut{$\phi$} }(\omega)\widehat{\ut{$h$}}(\omega)$ has a continuous and bounded inverse Fourier. 
Next, we show that the inverse Fourier of $\omega^{\beta}\indc_{\omega>1}$ is given by
\begin{align}\label{eq:InvFourier0}
&\mathcal{F}^{-1}\Big\{ \omega^{\beta} \indc_{\omega>1} \Big\}(x) = \tfrac{1}{2\pi} \frac{\Gamma(\beta+1 , -\jj x)}{(-\jj x)^{\beta+1}} \nonumber\\
&~~~~~ = 
\left\{\begin{array}{ll}
\tfrac{1}{2\pi} \frac{\Gamma(\beta+1)}{ (-\jj x)^{\beta+1} } -\tfrac{1}{2\pi} \sum_{m=0}^{\infty} \frac{ (-1)^m (- \jj x)^m}{ m! (\beta+1+m)}, & \beta \in \C\setminus \Z^{-},\\
&\\
\tfrac{-1}{2\pi}\big(\gamma + \ln |x| - \jj\tfrac{\pi}{2} {\rm sign}(x) + \sum_{m=1}^{\infty}\frac{(-\jj x)^m}{m\, m!}\big), & \beta= -1,
\end{array}
\right.
\end{align}
where $\Gamma(\cdot,\cdot)$ is the upper branch of the incomplete gamma function, and $\gamma$ is the Euler-Mascheroni constant. We first consider the case $\beta\not\in \Z^{-}$. Since $\omega^{\beta}\indc_{\omega>1}$ is not necessarily in $L_1$, we need to apply some techniques common for deriving the Fourier transform of tempered generalized function: 
\begin{align}\label{eq:InvFourier1}
\mathcal{F}^{-1}\Big\{ \omega^{\beta} \indc_{\omega>1}\Big\}(x) 
&= \tfrac{1}{2\pi} \int_{1}^{\infty} \omega^{\beta} \ee^{\jj \omega x} \dd \omega 
= \tfrac{1}{2\pi}\lim_{\mu \to 0^{+}} \int_{1}^{\infty} \omega^{\beta} \ee^{-(\mu-\jj x)\omega} \dd \omega \nonumber\\
&= \tfrac{1}{2\pi}\lim_{\mu \to 0^{+}} \lim_{M\to +\infty }\int_{1}^{M} \omega^{\beta} \ee^{-(\mu-\jj x)\omega} \dd \omega \nonumber\\
&= \tfrac{1}{2\pi}\lim_{\mu \to 0^{+}} \lim_{M\to +\infty } \frac{1}{(\mu-\jj x)^{\beta+1}} \int_{\CC} \; \tau^{\beta} \ee^{-\tau} \dd \tau, 
\end{align}
where $\CC$ stands for the finite line that connects the two points $\mu-\jj x$ and $(\mu-\jj x)M$ in the complex plane, and the latter integral is interpreted as a contour integration. Note that $\CC$ lies strictly on the right side of the imaginary axis, and  $z^{\beta} \ee^{-z}$ is analytic in this region. Further, $\gamma(\beta+1,z)$ (the lower branch of the incomplete gamma function) is an anti-derivative for this function ($\beta\not\in \Z^{-}$). Thus, \eqref{eq:InvFourier1} can be rewritten as
\begin{align}\label{eq:InvFourier2}
\mathcal{F}^{-1}\Big\{ \omega^{\beta} \indc_{\omega>1}\Big\}(x) 
&= \tfrac{1}{2\pi}\lim_{\mu \to 0^{+}} \lim_{M\to +\infty } \frac{\gamma\big(\beta+1 , (\mu-\jj x)M\big) - \gamma\big(\beta+1 , \mu-\jj x\big) }{(\mu-\jj x)^{\beta+1}}  .
\end{align}
It is known that $\lim_{|\zeta|\to\infty} \gamma(s,\zeta)=\Gamma(s)$, given that $|\measuredangle \zeta|<\tfrac{\pi}{2}$ while $|\zeta|\to\infty$. Indeed, this is the case for $\zeta=(\mu-\jj x)M$ as $M\to+\infty$. Therefore,
\begin{align}\label{eq:InvFourier3}
\mathcal{F}^{-1}\Big\{ \omega^{\beta} \indc_{\omega>1}\Big\}(x) 
&= \tfrac{1}{2\pi}\lim_{\mu \to 0^{+}}  \frac{\Gamma(\beta+1) - \gamma\big(\beta+1 , \mu-\jj x\big) }{(\mu-\jj x)^{\beta+1}}  \nonumber\\
&= \tfrac{1}{2\pi} \frac{\Gamma(\beta+1) - \gamma(\beta+1 , -\jj x) }{(-\jj x)^{\beta+1}} \nonumber\\
&= \tfrac{1}{2\pi} \frac{ \Gamma(\beta+1)  }{(-\jj x)^{\beta+1}} - \tfrac{1}{2\pi} \sum_{m=0}^{\infty} \frac{(-1)^m (-\jj x)^{m}}{ m! (\beta+1+m) } ,
\end{align}
which proves the first part of the claim in \eqref{eq:InvFourier0}. Note that $\beta=-1$ does not satisfy many of the results used in the above argument. Hence, we consider this case separately:
\begin{align*}
\mathcal{F}^{-1}\Big\{ \tfrac{ \indc_{\omega>1}}{\omega}\Big\}(x) = \tfrac{1}{2\pi} \int_{1}^{\infty} \frac{\ee^{\jj \omega t}}{\omega} \dd \omega =  \tfrac{1}{2\pi} E_1(-\jj t),
\end{align*}
where $E_1$ is the analytic continuation of the exponential integral function. The expression in \eqref{eq:InvFourier0} is in fact a well-known series expansion of $E_1(\zeta)$ at $\zeta = -\jj t$.

The inverse Fourier transform of  $|\omega|^{\beta} \indc_{\omega<-1}$ could also be obtained from \eqref{eq:InvFourier0} via
\begin{align}\label{eq:InvFourier4}
\mathcal{F}^{-1}\Big\{ |\omega|^{\beta} \indc_{\omega<-1}\Big\}(x) 
&= \mathcal{F}^{-1}\Big\{ \omega^{\beta} \indc_{\omega>1}\Big\}(-x) ,
\end{align}
which is valid due to the axis-flipping property of the Fourier transform ($\omega \rightarrow -\omega$).

We now use \eqref{eq:InvFourier0} to characterize the inverse Fourier of  $\widehat{\varphi}_{a,b}^{\Dord; \scalebox{0.6}{$k$}}(\omega)$ in \eqref{eq:PhiH}. Except the term $\widehat{ \ut{$\phi$} }(\omega) \widehat{\ut{$h$}}(\omega) $ which has a bounded and continuous inverse Fourier, we have that
\begin{align}\label{eq:ResFourier1}
\mathcal{F}^{-1}\Big\{& \tfrac{ \widehat{\varphi}^{(k-1)}(0) }{ (k-1) ! } |\omega|^{\lambda-1+\jj\imag(\Dord)} \big(\tfrac{\indc_{\omega>1}}{a} + \tfrac{(-1)^{k-1}\indc_{\omega<-1} }{b} \big) \Big\}(x)
\nonumber\\
& = 
\tfrac{ \widehat{\varphi}^{(k-1)}(0) }{ 2\pi \, (k-1) ! } 
\bigg( \underbrace{\frac{\Gamma(\lambda+\jj \imag(\Dord) \,,\, -\jj x)}{a\, (-\jj x)^{\lambda+\jj\imag (\Dord)}}
+  \frac{(-1)^{k-1} \Gamma(\lambda+\jj \imag(\Dord) \,,\, \jj x)}{b\,(\jj x)^{\lambda+\jj\imag(\Dord)}} }_{\Xi} \bigg),
\end{align}
where for $\lambda\neq 0$ or $\imag(\Dord) \neq 0$ (equivalently $\ut{$\Dord$}\neq 0$), $\Xi$ can be expressed as
\begin{align}\label{eq:ResFourier2}
\tfrac{\Gamma(\lambda+\jj\imag(\Dord)) }{a(-\jj x)^{\lambda+\jj\imag(\Dord)}} + \tfrac{(-1)^{k-1}\Gamma(\lambda+\jj\imag(\Dord)) }{b(\jj x)^{\lambda+\jj\imag(\Dord)}}
- \sum_{m=0}^{\infty} (-1)^m\tfrac{b (-\jj x)^m + (-1)^{k-1}a (\jj x)^m}{a\,b\,m!\,(\lambda+\jj\imag(\Dord) + m)},
\end{align}
and for $\lambda=\imag(\Dord)=0$  (equivalently, $\ut{$\Dord$}= 0$) as
\begin{align}\label{eq:ResFourier3}
-\tfrac{(b+a(-1)^{k-1})(\gamma + \ln |x|)}{a\,b} 
+ \tfrac{\jj \pi {\rm sign}(x) (b- (-1)^{k-1}a)}{2ab}
- \sum_{m=1}^{\infty} (-1)^m\tfrac{b (-\jj x)^m + a(-1)^{k-1} (\jj x)^m}{a\,b\,m\; m!}.
\end{align}
Because of the fact that
\begin{align*}
\bigg| \sum_{m=1}^{\infty} (-1)^m\tfrac{b (-\jj x)^m + (-1)^{k-1}a (\jj x)^m}{a\,b\,m!\,(\lambda+\jj\imag(\Dord) + m)} \bigg| 
\leq \tfrac{|a| + | b|}{ |a\,b\,( \lambda+1 + \jj \imag(\Dord) )| } \underbrace{\sum_{m=1}^{\infty} \tfrac{|x|^m}{m!} }_{\ee^{|x|}-1},
\end{align*}
we conclude that $\sum_{m=0}^{\infty} (-1)^m\tfrac{b (-\jj x)^m + (-1)^{k-1}a (\jj x)^m}{a\,b\,m!\,(\lambda+\jj\imag(\Dord) + m)}$ converges to a continuous function. 
Further, since $\pm\jj x = \exp\big(\log |x| \pm \jj \tfrac{\pi}{2} {\rm sign}(x) \big)$, we know that
\begin{align}\label{eq:IncGammaDenom}
\tfrac{1}{(\pm \jj x)^{\lambda+\jj \imag(\Dord)}} = 
\exp \hspace{-1mm}\Big( \hspace{-1.5mm} -\hspace{-1mm} \lambda \log |x| \hspace{-0.5mm} \pm \hspace{-0.5mm} \tfrac{\pi}{2} \imag(\Dord) \, {\rm sign}(x) \hspace{-1mm} - \hspace{-1mm} \jj \big( \hspace{-1mm} \imag(\Dord) \, \log |x|  \hspace{-0.5mm} \pm \hspace{-0.5mm} \tfrac{\pi}{2}\lambda \, {\rm sign}(x) \big) \hspace{-1mm} \Big)
\end{align}
is also a continuous function except possibly at $x=0$. 
Overall, we conclude that ${\varphi}_{a,b}^{\Dord; \scalebox{0.6}{$k$}}(x)$ in \eqref{eq:PhiH} (the inverse Fourier of $\widehat{\varphi}_{a,b}^{\Dord; \scalebox{0.6}{$k$}}$) is well-defined as a function and is continuous everywhere except possibly at $x= 0$; 
around $x= 0$, the behavior of ${\varphi}_{a,b}^{\Dord; \scalebox{0.6}{$k$}}(x)$ is determined by the $\frac{1}{(\pm \jj x)^{\lambda+\jj \imag(\Dord)}}$ terms for $\ut{$\Dord$}\neq 0$, and $-\tfrac{b+a(-1)^{k-1}}{a\,b} \ln |x|
+ \tfrac{\jj \pi  (b- (-1)^{k-1}a)}{2ab} {\rm sign}(x) $ term for $\ut{$\Dord$}=0$.
To proceed, we check different cases separately.

\begin{enumerate}[(i)]
\item $\real(\Dord)\not\in \Z$ and $k= \lfloor \real(\Dord)\rfloor $. Thus, $-1 < \lambda <0$. In this case, 
\begin{align*}
\lim_{|x|\to 0}\left|\frac{1}{(\pm \jj x)^{\lambda+\jj \imag(\Dord)}}\right| =0,
\end{align*}
which results in 
$$\lim_{x\rightarrow 0} \frac{\Gamma(\lambda+\jj \imag(\Dord),\pm\jj x)}{(\pm\jj x)^{\lambda+\jj\imag(\Dord)}} = \frac{-1}{\lambda+\jj \imag(\Dord)}.$$
Consequently, 
$\tfrac{\Gamma(\lambda+\jj \imag(\Dord), \pm\jj x)}{(\pm\jj x)^{\lambda+\jj\imag (\Dord)}}$ is continuous everywhere including at $x=0$. This proves claim (\ref{claim:k_1_NonZ}).

\item $\real(\Dord)\in \Z$, $\imag(\Dord)\neq 0$, and $k= \lfloor \real(\Dord)\rfloor $. Thus, $\lambda=0$ but $\ut{$\Dord$}\neq 0$. Recalling \eqref{eq:IncGammaDenom}, we observe that
\begin{align*}
\bigg|\frac{1}{(\jj x)^{\jj\imag(\Dord)}}\bigg| 
 \leq \ee^{\tfrac{\pi}{2}|\imag(\Dord)|}.
\end{align*}
Consequently, $\tfrac{\Gamma(\jj \imag(\Dord),\jj x)}{(\jj x)^{\jj\imag(\Dord)}}$ (and in turn, ${\varphi}_{a,b}^{\Dord; \scalebox{0.6}{$k$}}$) is bounded around $x=0$. However, the oscillating nature of $\tfrac{\Gamma(\jj \imag(\Dord),\jj x)}{(\jj x)^{\jj\imag(\Dord)}}$ around $x=0$ makes it discontinuous.

\item $\Dord\in \Z$, $h_{a,b}^{\Dord}(\omega) \equiv c\,\omega^{\Dord}$, and $k= \lfloor \real(\Dord)\rfloor $. This implies that $\ut{$\Dord$}= 0$ and $b=c\, (-1)^k= (-1)^k \, a$. According to \eqref{eq:ResFourier3}, the term $\ln|x|$ vanishes. Therefore, the inverse Fourier remains bounded around $x=0$; however, due to the existence of the ${\rm sign}(x)$ term, it will be discontinuous at $x=0$.

\item $\Dord\in \Z$, $h_{a,b}^{\Dord}(\omega) \not\equiv c\,\omega^{\Dord}$, and $k= \lfloor \real(\Dord)\rfloor $. This case is very similar to the previous case, except that $b \neq (-1)^k \, a$ and the $\log|x|$ term remains. Thus, ${\varphi}_{a,b}^{\Dord; \scalebox{0.6}{$k$}}(x)$ is proportional to $\log|x|$ around $x=0$.

\item $k> \lfloor \real(\Dord)\rfloor$. As a result $0<\lambda=k-\real(\Dord) $ and obviously $\ut{$\Dord$}\neq 0$. This implies that $\tfrac{1}{(\mp\jj x)^{\lambda+\jj \imag(\Dord)}}$, and as a result ${\varphi}_{a,b}^{\Dord; \scalebox{0.6}{$k$}}(x)$, are singular at the origin such that $| x|^{\lambda}|{\varphi}_{a,b}^{\Dord; \scalebox{0.6}{$k$}}(x)|$ is bounded and discontinuous around $x=0$. Note that, there is no choice of $a,b$ such that the two singularities in \eqref{eq:ResFourier2} completely cancel each other out (the cancellation can happen either for $x>0$ or for $x<0$). 
$\hspace{\stretch{1}}\blacksquare$
\end{enumerate}

\begin{lemma}\label{lemma:InverseFourier}
Let $\Dord\in\C\setminus\Z$,
then, for all $k\geq 1$ we have that
\begin{align}\label{eq:InverseFourier}
2\pi \mathcal{F}^{-1}\bigg\{ \omega_{+}^{\Dord-k}\Big(\ee^{\jj\omega t} - \sum_{j=0}^{k-1} \tfrac{(\jj \omega t)^j}{j!}\Big) \bigg\}(x) 
= \tfrac{\Gamma(\Dord-k+1)}{\big(-\jj(x+t)\big)^{\Dord-k+1}} -\sum_{j=0}^{k-1} \tfrac{(\jj t)^{k-j-1} \Gamma(\Dord-j)}{(k-j-1)! (-\jj x)^{\Dord-j}}.
\end{align}
\end{lemma}

\textbf{Proof.}
For non-integer $\Dord$, it is know that \cite{JonesBook1982,Forster2006}
\begin{align}\label{eq:InvF_Wp}
2\pi\mathcal{F}^{-1}\big\{ w_{+}^{\Dord} \big\}(x) =\frac{\Gamma(\Dord+1)}{(-\jj x)^{\Dord+1}} .
\end{align}
The above result is the key to prove the claim in Lemma \ref{lemma:InverseFourier}:
\begin{align*}
&\tfrac{\Gamma(\Dord-k+1)}{\big(-\jj(x+t)\big)^{\Dord-k+1}} -\sum_{j=0}^{k-1} \tfrac{(\jj t)^{k-j-1} \Gamma(\Dord-j)}{(k-j-1)! (-\jj x)^{\Dord-j}} \nonumber\\
&~~~ = 2\pi\mathcal{F}^{-1}\big\{ w_{+}^{\Dord-k} \big\}(x+t) - 2\pi\sum_{j=0}^{k-1}
\tfrac{(\jj t)^{k-j-1}}{(k-j-1)!}  \mathcal{F}^{-1}\big\{ w_{+}^{\Dord-j-1} \big\}(x) \nonumber\\
&~~~ = 2\pi\mathcal{F}^{-1}\big\{ w_{+}^{\Dord-k} \ee^{\jj \omega t}\big\}(x) - 2\pi\sum_{j=0}^{k-1}  \mathcal{F}^{-1}\Big\{\tfrac{(\jj t)^{k-j-1}}{(k-j-1)!} w_{+}^{\Dord-j-1} \Big\}(x)\nonumber\\
&~~~ = 2\pi \mathcal{F}^{-1}\Big\{ w_{+}^{\Dord-k} \ee^{\jj \omega t} - \sum_{j=0}^{k-1} \tfrac{(\jj t)^{k-j-1}}{(k-j-1)!} w_{+}^{\Dord-j-1}\Big\}(x)\nonumber\\
&~~~ = 2\pi \mathcal{F}^{-1}\bigg\{ w_{+}^{\Dord-k} \Big(\ee^{\jj \omega t} - \sum_{j=0}^{k-1} \tfrac{(\jj \omega t)^{k-j-1}}{(k-j-1)!} \Big)\bigg\}(x),
\end{align*}
which is the same as claim. $\hspace{\stretch{1}}\blacksquare$


\section{Proofs}\label{sec:Proofs}

In this section, we prove the results stated in Section \ref{sec:Theos}. In the proofs we make use of the lemmas in Section \ref{sec:Lemmas}.

\subsection{Proof of Theorem \ref{theo:decayRate}} 
As $\Dop_{a,b}^{\Dord}$ corresponds to a Fourier multiplier, it is a convolutional operator. Since $\varphi\in\SS$ is infinitely differentiable, this also carries over to the convolution of $\varphi$ with any generalized function. Thus, $\big(\Dop_{a,b}^{\Dord}\varphi\big)(x)$ is infinitely differentiable.

The statement (\ref{decayTheo_claim:schwarz}) is a classical result which simply follows from the definition of Schwartz functions.
To prove other statements, let $T\in[2,4[$ and consider the following input to the system:
$$\varphi(x)-T\,\varphi(Tx) ~ \stackrel{\Dop_{a,b}^{\Dord}}{\longmapsto} ~ \varphi_{a,b}^{\Dord}(x) - T^{\Dord+1} \varphi_{a,b}^{\Dord}(Tx).$$
Using the Fourier representation we have that
\begin{align}\label{eq:SingleDiff}
\varphi_{a,b}^{\Dord}(x) - T^{\Dord+1} \varphi_{a,b}^{\Dord}(Tx) &= \tfrac{1}{2\pi}\int_{\R} \big( \widehat{\varphi}(\omega) - \widehat{\varphi}(\tfrac{\omega}{T}) \big) h_{a,b}^{\Dord}(\omega) \ee^{\jj\omega x}\dd \omega.
\end{align}
Since $\varphi\in\mathcal{S}$, it is well-known that $\widehat{\varphi}\in\mathcal{S}$. Hence, $\widehat{\varphi}$ is infinitely diffefrentiable and $\big|\widehat{\varphi}^{(j)}(\omega)\big|\big(1+|\omega|^{r+j}\big)$ is bounded for all $1\leq j,r$. 
According to Lemma \ref{lemma:L1property}, $\big( \widehat{\varphi}(\omega) - \widehat{\varphi}(\tfrac{\omega}{T}) \big)h_{a,b}^{\Dord}(\omega)$ is $n_{\Dord}=\lceil \real(\Dord)\rceil + 1$ times continuously differentiable and $\frac{\dd^m}{\dd \omega^m} \Big(\big( \widehat{\varphi}(\omega) - \widehat{\varphi}(\tfrac{\omega}{T}) \big)h_{a,b}^{\Dord}(\omega)\Big)$ has finite $L_1$ norm for all $1\leq m\leq n_{\Dord}$ and $T\in[2,4[$\footnote{Note that when $h_{a,b}^{\Dord}(\omega)=c \omega^n$ for some $n\in\N$ and $c\in\C$, $\big( \widehat{\varphi}(\omega) - \widehat{\varphi}(\tfrac{\omega}{T}) \big)h_{a,b}^{\Dord}(\omega)$ is infinitely differentiable, and these are the only cases with this property}. 
Hence, by integration by parts, we can rewrite \eqref{eq:SingleDiff} as
\begin{align}\label{eq:SingleDiff2}
\varphi_{a,b}^{\Dord}(x) - T^{\Dord+1} \varphi_{a,b}^{\Dord}(Tx) &= \tfrac{1}{2\pi} (\tfrac{- 1}{x})^{n_{\Dord}} \int_{\R} \underbrace{\tfrac{\dd^{n_{\Dord}}}{\dd \omega^{n_{\Dord}}}\Big( \big( \widehat{\varphi}(\omega) - \widehat{\varphi}(\tfrac{\omega}{T}) \big)h_{a,b}^{\Dord}(\omega) \Big)}_{\widehat{\psi}_{_{T}}(\omega) } \ee^{\jj\omega x}\dd \omega \nonumber\\
&= \tfrac{1}{2\pi} (\tfrac{- 1}{x})^{n_{\Dord}} \int_{\R} \widehat{\psi}_{_{T}}(\omega) \ee^{\jj\omega x}\dd \omega = (\tfrac{- 1}{x})^{n_{\Dord}} \psi_{_{T}}(x).
\end{align}
Again based on Lemma \ref{lemma:L1property}, we have that $\|\widehat{\psi}_{_{T}}\|_1\leq c 4^{\real(\Dord) + 1 -\lceil \real (\Dord)\rceil}$. Thus,
\begin{align*}
\big|\varphi_{a,b}^{\Dord}(x) - T^{\Dord+1} \varphi_{a,b}^{\Dord}(Tx)\big| \leq \frac{\|\psi_{_{T}}\|_{\infty}}{|x|^{n_{\Dord}}} =  \frac{ \|\widehat{\psi}_{_{T}}\|_1 }{|x|^{n_{\Dord}}} 
<\frac{c\; 4^{\real(\Dord) +1 -\lceil \real(\Dord)\rceil}}{|x|^{n_{\Dord}}}.
\end{align*}
This implies that $\varphi_{a,b}^{\Dord}$ satisfies the constraint of Lemma \ref{lemma:ArashTechnique} with $\Dord$ and $m$ in Lemma  \ref{lemma:ArashTechnique} being replaced with $\Dord+1$ and $n_{\Dord}$, respectively. Thus,  $\bar{c}\in\R^{+}$ exists such that
\begin{align}\label{eq:TailB1}
\forall \, x\in \R\setminus\{0\}:~~ \big|\varphi_{a,b}^{\Dord}(x)\big| \leq \left\{ \begin{array}{cl}
\frac{\bar{c}}{ |x|^{\real\,(\Dord+1)} }, & n_{\Dord}\neq \real\,(\Dord+1),\phantom{\Big|}\\
\bar{c} \, \frac{1+\big|\log|x|\big|}{ |x|^{\real\,(\Dord+1)} }, & n_{\Dord}= \real\,(\Dord+1), \phantom{\Big|}
\end{array}
\right.
\end{align}
where we used the fact that $\min(|x|^{\real(\Dord+1)}\,,\, |x|^{n_{\Dord}}) = |x|^{\real(\Dord+1)}$, as $n_{\Dord}=\lceil \real(\Dord)\rceil + 1\geq \real (\Dord+1)$. We also recall that $\real (\Dord)>-1$; thus, $\widehat{\varphi}(\omega) h_{a,b}^{\Dord}(\omega)$ is locally integrable. Further, $\widehat{\varphi}\in\SS$ has rapid decay as $|\omega|\to\infty$. These two properties imply that $\widehat{\varphi}(\omega) h_{a,b}^{\Dord}(\omega)\in L_1$. Consequently, $\varphi_{a,b}^{\Dord}(x)$ (the inverse Fourier of $\widehat{\varphi}(\omega)h_{a,b}^{\Dord}(\omega)$) is both bounded and continuous:
\begin{align}\label{eq:TailB2}
\exists\, M_{\varphi}\in\R^{+},~ \forall\,x\in\R:~~~ \big|\varphi_{a,b}^{\Dord}(x)\big| \leq M_{\varphi}.
\end{align}

For $\real(\Dord)\not\in\Z$, we have that $n_{\Dord}\neq \real\,(\Dord+1)$, and can combine \eqref{eq:TailB1} and \eqref{eq:TailB2} as 
\begin{align}\label{eq:TailB3}
\big|\varphi_{a,b}^{\Dord}(x)\big| 
\leq 
\min\big( \tfrac{\bar{c}}{ |x|^{\real\,(\Dord+1)} } \,,\, M_{\varphi} \big) 
&\leq 
\tfrac{2}{ \tfrac{ |x|^{\real\,(\Dord+1)} }{\bar{c}} + \tfrac{1}{M_{\varphi}} }  
= 
\tfrac{2 \bar{c}}{ |x|^{\real\,(\Dord+1)} + \tfrac{\bar{c}}{M_{\varphi}}} \nonumber\\
&\leq 
\tfrac{2 (\bar{c} + M_{\varphi})}{ 1 + |x|^{\real\,(\Dord+1)} },
\end{align}
which is the same bound as claimed in statement (\ref{decayTheo_claim:gooddecay}). For $\real(\Dord)\in\Z$ which coincides with $n_{\Dord}= \real\,(\Dord+1)$, we have that:
\begin{align}\label{eq:TailB4}
\big|\varphi_{a,b}^{\Dord}(x)\big| 
\leq 
\min\Big( \bar{c} \, \tfrac{1+\big|\log|x|\big|}{ |x|^{\real\,(\Dord+1)} } \,,\, M_{\varphi} \Big) 
\leq 
\tfrac{2 (\bar{c} + M_{\varphi})\log(1+|x|)}{ 1 + |x|^{\real\,(\Dord+1)} },
\end{align}
which is again the required bound in statement (\ref{decayTheo_claim:slowdecay}). The validity of \eqref{eq:TailB4} could be verified via
\begin{align*}
\tfrac{2 (\bar{c} + M_{\varphi})\log(1+|x|)}{ 1 + |x|^{\real\,(\Dord+1)} } \geq \left\{\begin{array}{cl}
M_{\varphi}, & |x|\leq 1, \phantom{\Big|}\\
\bar{c} \, \tfrac{1+\big|\log|x|\big|}{ |x|^{\real\,(\Dord+1)} }, & |x| \geq 1. \phantom{\Big|}
\end{array}
\right.
\end{align*}
Based on \eqref{eq:TailB3} and \eqref{eq:TailB4}, it is now easy to show that $\varphi_{a,b}^{\Dord} \in L_p$ for $p>\frac{1}{\real(\Dord) + 1}$. 
The continuity of the mapping from $\SS$ to $L_p$ also follows from linearity and boundedness of the mapping.
$\hspace{\stretch{1}}\blacksquare$


\subsection{Proof of Proposition \ref{prop:Integration_Adjoint}}  
We first show that $\Iop_{a,b}^{\Dord; \scalebox{0.6}{$k$}}$ is well-defined for $\varphi\in\mathcal{S}$. For this purpose, we express $\big( \Iop_{a,b}^{\Dord; \scalebox{0.6}{$k$}}\varphi\big)(x)$ as
\begin{align}\label{eq:nonLSI_integral}
\big( \Iop_{a,b}^{\Dord; \scalebox{0.6}{$k$}}\varphi\big)(x)
=& 
 \tfrac{1}{2\pi}\int_{\R} \widehat{\varphi}( \omega) 
\Big(\ee^{\jj \omega x}-\sum_{j=0}^{k-1}\tfrac{(\jj x)^j}{j!} \omega^j  \Big)\tfrac{1}{h_{a,b}^{\Dord}(\omega)}\, \dd \omega \nonumber\\
= &
 \tfrac{1}{2\pi}\int_{\R} \widehat{\varphi}( \omega) 
\underbrace{ \frac{\ee^{\jj \omega x}-\sum_{j=0}^{k-1}\tfrac{(\jj x)^j}{j!} \omega^j  }{w^{k} }}_{\widehat{\phi}(\omega)}
\underbrace{ \tfrac{\omega^{k}}{h_{a,b}^{\Dord}(\omega)} }_{ {\ut{$h$}}(\omega) } \, \dd \omega.
\end{align}
Lemma \ref{lemma:ContDeriv} implies that $\widehat{\phi}(\omega)$ is continuous and decays no slower than $\frac{1}{1+|\omega|}$. 
As $\varphi\in\SS$, $\widehat{\varphi}$ is also in $\SS$. Specially, $\widehat{\varphi}$ is continuous everywhere and asymptotically decays faster than $\frac{1}{1+|\omega|^{r}}$ for any $r>0$. 
It is easy to check that ${\ut{$h$}}(\omega)=h_{\scalebox{0.6}{\ut{$a$}},\scalebox{0.6}{\ut{$b$}}}^{\scalebox{0.6}{\ut{$\Dord$}}}(\omega)$ is homogeneous of degree $\ut{$\Dord$} = k-\Dord$ with $-1<\real(\ut{$\Dord$})$, and $\ut{$a$} = \frac{1}{a},~ \ut{$b$}=\frac{(-1)^{k}}{b}$. 
Hence, $\big|\widehat{\varphi}(\omega) \widehat{\phi}(\omega)\widehat{\ut{$h$}}(\omega)\big|$ could be upper-bounded by $c |\omega|^{\real(\ut{$\Dord$})}$ around $\omega=0$, and by $c |\omega|^{ - 2}$ as $|\omega|\to \infty$ ($\widehat{\varphi}$ has a super-polynomial decay rate). Consequently, $\big|\widehat{\varphi}(\omega) \widehat{\phi}(\omega)\widehat{\ut{$h$}}(\omega)\big|$ could be fully dominated by $c \min\big( |\omega|^{\real(\ut{$\Dord$})} \,,\, |\omega|^{-2}\big)$ when $c\in\R^+$ is large enough. The latter has a finite integral over the real line, which shows that the integral in \eqref{eq:nonLSI_integral} is well-defined. 
To see that $\Iop_{a,b}^{\Dord; \scalebox{0.6}{$k$}}$ is scale invariant, pick any $T>0$ and consider 
\begin{align*}
\varphi (Tx)
\stackrel{\Iop_{a,b}^{\Dord; \scalebox{0.6}{$k$}}}{\longmapsto} 
& \phantom{=}
 \tfrac{1}{2\pi}\int_{\R} \tfrac{1}{T}\widehat{\varphi}\big( \tfrac{\omega}{T}\big) 
\tfrac{\ee^{\jj \omega x}- \big(\sum_{j=0}^{k-1}\tfrac{(\jj x)^j}{j!} \omega^j \big) }{h_{a,b}^{\Dord}(\omega)} \, \dd \omega \nonumber\\
&=
 \tfrac{1}{2\pi}\int_{\R} \widehat{\varphi}( \zeta) 
\tfrac{\ee^{\jj \zeta T x} - \big(\sum_{j=0}^{k-1}\tfrac{(\jj x)^j}{j!} (T\zeta)^j \big) }{h_{a,b}^{\Dord}(T \zeta)} \, \dd \zeta 
\nonumber\\
&=
 \tfrac{T^{-\Dord}}{2\pi}\int_{\R} \widehat{\varphi}( \zeta) 
\tfrac{\ee^{\jj \zeta T x} - \big(\sum_{j=0}^{k-1}\tfrac{(\jj T x)^j}{j!} \zeta^j \big) }{h_{a,b}^{\Dord}(T \zeta)} \, \dd \zeta = T^{-\Dord} \big(\Iop_{a,b}^{\Dord; \scalebox{0.6}{$k$}} \varphi\big)( Tx).
\end{align*}

Next, we show that $\Iop_{a,b}^{\Dord; \scalebox{0.6}{$k$}}$ is the right-inverse 
of the LSI operator $\Dop_{a,b}^{\Dord}$ when $k=\lfloor \real(\Dord)\rfloor$ or $k=\lceil \real(\Dord)\rceil$:
\begin{align}\label{eq:RightInverse}
\big(\Dop_{a,b}^{\Dord} \, \Iop_{a,b}^{\Dord; \scalebox{0.6}{$k$}} \, \varphi \big)(x) 
=&
 \mathcal{F}^{-1}_{\omega}\bigg\{h_{a,b}^{\Dord}(\omega) \mathcal{F}_{\tau}\Big\{ \big( \Iop_{a,b}^{\Dord; \scalebox{0.6}{$k$}} \, \varphi \big)(\tau) \Big\}(\omega)\bigg\}(x) \nonumber\\
 =&
 \mathcal{F}^{-1}_{\omega}\bigg\{ h_{a,b}^{\Dord}(\omega)   \mathcal{F}_{\tau}\Big\{ \tfrac{1}{2\pi} \int_{\R} \widehat{\varphi}(\zeta) \tfrac{\ee^{\jj \zeta \tau} - \sum_{j=0}^{k-1}\tfrac{(\jj  \tau)^j}{j!} \zeta^j}{h_{a,b}^{\Dord}(\zeta)}
  \dd \zeta \Big\}(\omega) \bigg\}(x) \nonumber\\
=&
 \tfrac{1}{2\pi} \mathcal{F}^{-1}_{\omega}\bigg\{ h_{a,b}^{\Dord}(\omega) \int_{\R} \tfrac{\widehat{\varphi}(\zeta)}{ h_{a,b}^{\Dord}(\zeta)}
 \mathcal{F}_{\tau}\Big\{\ee^{\jj \zeta \tau} - \sum_{j=0}^{k-1}\tfrac{(\jj  \tau)^j}{j!} \zeta^j
  \Big\}(\omega) \dd \zeta \bigg\}(x) \nonumber\\
=&
  \mathcal{F}^{-1}_{\omega}\bigg\{ h_{a,b}^{\Dord}(\omega)    \int_{\R} \tfrac{\widehat{\varphi}(\zeta)}{ h_{a,b}^{\Dord}(\zeta)}
 \Big(\delta(\omega-\zeta) - \sum_{j=0}^{k-1}\tfrac{1}{j!} \zeta^j \delta^{(j)}(\omega)
  \Big) \dd \zeta \bigg\}(x) \nonumber\\
=
  \mathcal{F}^{-1}_{\omega}\bigg\{   & \int_{\R} \tfrac{\widehat{\varphi}(\zeta)}{ h_{a,b}^{\Dord}(\zeta)}
 \bigg(h_{a,b}^{\Dord}(\omega) \delta(\omega-\zeta)   \hspace{-1mm}-\hspace{-1mm} \sum_{j=0}^{k-1}\tfrac{1}{j!} \zeta^j h_{a,b}^{\Dord}(\omega) \delta^{(j)}(\omega) 
  \bigg) \dd \zeta \bigg\}(x).
\end{align}
We know that $\frac{\dd^{j}}{\dd \omega^j} h_{a,b}^{\Dord}(\omega)\Big|_{\omega=0} = 0$ for all $0\leq j < \real(\Dord)$; this means that $ h_{a,b}^{\Dord}(\omega) \delta^{(j)}(\omega)\equiv 0$ for $0\leq j < \real(\Dord)$. In case $k=\lfloor \real(\Dord)\rfloor$ or $k=\lceil \real(\Dord)\rceil$, the interval $\in[0 , \real(\Dord)[$ covers the full range of required $j$ values in the summation in \eqref{eq:RightInverse}. Hence, this summation could be dropped from \eqref{eq:RightInverse}:
\begin{align}\label{eq:RightInverse2}
\big(\Dop_{a,b}^{\Dord} \, \Iop_{a,b}^{\Dord; \scalebox{0.6}{$k$}} \, \varphi \big)(x) 
&=
  \mathcal{F}^{-1}_{\omega}\bigg\{   \int_{\R} \tfrac{\widehat{\varphi}(\zeta)}{ h_{a,b}^{\Dord}(\zeta)}
 h_{a,b}^{\Dord}(\omega) \delta(\omega-\zeta)  
 \dd \zeta \bigg\}(x) \nonumber\\
 & =
   \mathcal{F}^{-1}_{\omega}\bigg\{   \widehat{\varphi}(\omega) 
 \tfrac{ h_{a,b}^{\Dord}(\omega) }{ h_{a,b}^{\Dord}(\omega) }   \bigg\}(x) 
 =
 \varphi(x),
\end{align}
which proves that $\Iop_{a,b}^{\Dord; \scalebox{0.6}{$k$}}$ is the right-inverse of $\Dop_{a,b}^{\Dord}$.

Our last task is to find the adjoint operator of $\Iop_{a,b}^{\Dord; \scalebox{0.6}{$k$}}$. 
Let $\varphi,\psi\in\SS$. We can write
\begin{align}\label{eq:Adjoint}
\big\langle \big(\Iop_{a,b}^{\Dord; \scalebox{0.6}{$k$}}\,\varphi\big)(x)\,,\, \psi(x) \big\rangle 
=&
\int_{\R} \big(\Iop_{a,b}^{\Dord; \scalebox{0.6}{$k$}}\,\varphi\big)(x) \, \psi(x)  \dd x \nonumber\\
=&
 \tfrac{1}{2\pi} \int_{\R} \int_{\R} \widehat{\varphi}(\omega) \widehat{h}(\omega) \tfrac{\ee^{\jj \omega x} - \sum_{j=0}^{k-1}\tfrac{(\jj x)^j}{j!} \omega^j}{ h_{a,b}^{\Dord}(\omega) }  \psi(x)  \dd \omega \dd x \nonumber\\
 =&
 \tfrac{1}{2\pi}  \int_{\R}  \tfrac{\widehat{\varphi}(\omega)}{ h_{a,b}^{\Dord}(\omega) } \bigg(\int_{\R}\Big(\ee^{\jj \omega x} - \sum_{j=0}^{k-1}\tfrac{(\jj x)^j}{j!} \omega^j\Big)  \psi(x)  \dd x  \bigg) \dd \omega .
\end{align}
We use Parseval's theorem to obtain
\begin{align}\label{eq:Adjoint1}
\int_{\R}\Big(\ee^{\jj \omega x} - \sum_{j=0}^{k-1}\tfrac{(\jj x)^j}{j!} \omega^j\Big)  \psi(x)  \dd x
=&
\tfrac{1}{2\pi}\int_{\R} \widehat{\psi}(\zeta) \mathcal{F}_x\Big\{ \ee^{\jj \omega x} - \sum_{j=0}^{k-1}\tfrac{(\jj x)^j}{j!} \omega^j \Big\}(-\zeta) \dd\zeta \nonumber\\
=&
\int_{\R} \widehat{\psi}(\zeta) \Big( \delta(-\zeta-\omega) - \sum_{j=0}^{k-1}\tfrac{\omega^j}{j!}\delta^{(j)}(\zeta)  \Big) \dd\zeta \nonumber\\
=&
\widehat{\psi}(-\omega) - \sum_{j=0}^{k-1} \tfrac{(-\omega)^j}{j!} \widehat{\psi}^{(j)}(0).
\end{align}
By rewriting \eqref{eq:Adjoint} based on \eqref{eq:Adjoint1}, we arrive at
\begin{align}\label{eq:Adjoint2}
\big\langle \big(\Iop_{a,b}^{\Dord; \scalebox{0.6}{$k$}}\,\varphi\big)(x) \,,\, \psi(x) \big\rangle 
=&
  \tfrac{1}{2\pi} \int_{\R}  \tfrac{\widehat{\varphi}(\omega)}{ h_{a,b}^{\Dord}(\omega)} \Big( \widehat{\psi}(-\omega) - \sum_{j=0}^{k-1} \tfrac{(-\omega)^j}{j!} \widehat{\psi}^{(j)}(0)  \Big) \dd \omega \nonumber\\
  =&
 \int_{\R}  \varphi(x)  \mathcal{F}^{-1}_{\omega}\bigg\{  \tfrac{ \widehat{\psi}(-\omega) - \sum_{j=0}^{k-1} \tfrac{(-\omega)^j}{j!} \widehat{\psi}^{(j)}(0) }{ h_{a,b}^{\Dord}(\omega) } \bigg\}(-x) \dd x \nonumber\\
  =&
 \int_{\R}  \varphi(x)  \mathcal{F}^{-1}_{\omega}\bigg\{  \tfrac{ \widehat{\psi}(\omega) - \sum_{j=0}^{k-1} \tfrac{\omega^j}{j!} \widehat{\psi}^{(j)}(0) }{ h_{a,b}^{\Dord}(-\omega) } \bigg\}(x) \dd x \nonumber\\
   =&
   \big\langle \varphi(x) \,,\, \big(\Iop_{a,b}^{(\Dord; \scalebox{0.6}{$k$})*}\,\psi\big)(x) \big\rangle ,
\end{align}
where $\big(\Iop_{a,b}^{(\Dord; \scalebox{0.6}{$k$})*}\,\psi\big)(x) = \mathcal{F}^{-1}_{\omega}\bigg\{  \tfrac{ \widehat{\psi}(\omega) - \sum_{j=0}^{k-1} \tfrac{\omega^j}{j!} \widehat{\psi}^{(j)}(0) }{ h_{a,b}^{\Dord}(-\omega) } \bigg\}(x)$.
$\hspace{\stretch{1}}\blacksquare$


\subsection{Proof of Theorem \ref{theo:Integrator1}}
First note that Lemma \ref{lemma:BasicsIntegral} implies that for $\varphi\in\SS$, the output $\big(\Iop_{a,b}^{\Dord; \scalebox{0.4}{$k$}*}\varphi\big)(x)$ is well-defined and continuous at $x\neq 0$; further, this lemma describes the behavior of $\big(\Iop_{a,b}^{\Dord; \scalebox{0.4}{$k$}*}\varphi\big)(x)$ around $x=0$ the same way as claimed in Theorem \ref{theo:Integrator1}. Thus, to complete the proof of Theorem \ref{theo:Integrator1}, we need to investigate the decay properties of $\big(\Iop_{a,b}^{\Dord; \scalebox{0.4}{$k$}*}\varphi\big)(x)$ and the $L_p$ spaces to which this function belongs.

For $k=\lfloor\real(\Dord) \rfloor $ (parts (\ref{stat:Integrator1_Integer})-(\ref{stat:Integrator1_NonInteger}) of Theorem \ref{theo:Integrator1}), Lemma \ref{lemma:BasicsIntegral} shows that $\big|\big(\Iop_{a,b}^{\Dord; \scalebox{0.4}{$k$}*}\varphi\big)(x)\big|$ can be upper-bounded by $\nu+\kappa\big|\log |x|\big|$ at $|x|\leq 1$ for some $\nu,\kappa\in\R^{+}\cup\{0\}$. Thus,
\begin{align}\label{eq:SmallK0}
\int_{-1}^{1} \big|\big(\Iop_{a,b}^{\Dord; \scalebox{0.4}{$k$}*}\varphi\big)(x)\big|^p \dd x 
\leq  2\,\kappa^p \int_{0}^{1} \big( \tfrac{\nu}{\kappa} + |\log x|\big)^p \dd x
= 2\,\kappa^p \ee^{\frac{\nu}{\kappa}} \, \Gamma (p+1 \,,\, \tfrac{\nu}{\kappa} )
< \infty,
\end{align}
for any $p>0$. 

For $k>\lfloor\real(\Dord) \rfloor $ (parts (\ref{stat:Integrator2_NonInteger})-(\ref{stat:Integrator2_Integer}) of Theorem \ref{theo:Integrator1}), Lemma \ref{lemma:BasicsIntegral} implies that $|x|^{k-\real(\Dord)}\big|\big(\Iop_{a,b}^{\Dord; \scalebox{0.4}{$k$}*}\varphi\big)(x)\big|$ is bounded. Thus,
\begin{align}\label{eq:LargeK0}
\int_{-1}^{1} \big|\big(\Iop_{a,b}^{\Dord; \scalebox{0.4}{$k$}*}\varphi\big)(x)\big|^p \dd x 
\leq  2\int_{0}^{1}  \tfrac{c}{|x|^{p(k-\real(\Dord))}}\dd x
< \infty,
\end{align}
for  $p<\frac{1}{k-\real(\Dord)}$. 

We continue with a similar technique as in the proof of Lemma \ref{lemma:BasicsIntegral}:
\begin{align}
\mathcal{F}\big\{ \big(\Iop_{a,b}^{\Dord; \scalebox{0.4}{$k$}*}\varphi\big)(x)\big\}(\omega) = \underbrace{\frac{\widehat{\varphi}(\omega)-\sum_{j=0}^{k-1}\tfrac{\widehat{\varphi}^{(j)}(0)}{j!}\omega^j}{\omega^{k}} }_{\widehat{\phi}(\omega)} \underbrace{ \frac{\omega^{k}}{h_{a,b}^{\Dord}(-\omega)}}_{{\ut{$h$}}(\omega)} = \widehat{\phi}(\omega) {\ut{$h$}}(\omega),
\end{align}
where we know that ${\ut{$h$}}(\omega) = h_{a',b'}^{\scalebox{0.6}{\ut{$\Dord$}}}(\omega)$ is a homogeneous function of degree $\ut{$\Dord$}=k-\Dord$ with $-1<\real (\ut{$\Dord$})$ (and $a'= \frac{1}{b}, \, b'=\frac{(-1)^k}{a}$). 
Because of $\varphi\in\SS$, $\widehat{\varphi}$ is also in the Schwartz space and is infinitely differentiable; moreover, $\frac{\widehat{\varphi}^{(j)}(\omega)}{1+|\omega|^m}$ is bounded for all $0\leq j,m$. 
As a result, Lemma \ref{lemma:ContDeriv} implies that $\widehat{\phi}$ is infinitely differentiable, and $\widehat{\phi}^{(n)}(\omega)(1+|\omega|^{n+1})$ is bounded for all $n\geq 0$. 

Next, we consider the output of $\Iop_{a,b}^{\Dord; \scalebox{0.4}{$k$}*}$ to $T^{k+1}\varphi(Tx)$, where $T\in[2,4[$:
\begin{align*}
T^{k+1}\varphi(Tx) \phantom{=} \stackrel{\Iop_{a,b}^{\Dord; \scalebox{0.4}{$k$}*}}{\longmapsto}  \phantom{=\;} & \mathcal{F}^{-1}\bigg\{ T^{k} \tfrac{\widehat{\varphi}(\tfrac{\omega}{T})-\sum_{j=0}^{k-1}\tfrac{\widehat{\varphi}^{(j)}(0)}{T^j}\tfrac{\omega^j}{j!} }{h_{a,b}^{\Dord}(-\omega)} \bigg\}(x) \nonumber\\
=\; & \mathcal{F}^{-1}\bigg\{  \frac{\widehat{\varphi}(\tfrac{\omega}{T})-\sum_{j=0}^{k-1}\tfrac{\widehat{\varphi}^{(j)}(0)}{j!}(\tfrac{\omega}{T})^j }{(\tfrac{\omega}{T})^{k}} \tfrac{\omega^{k}}{h_{a,b}^{\Dord}(-\omega)} \bigg\}(x)\nonumber\\
=\; & \mathcal{F}^{-1}\Big\{  \widehat{\phi}(\tfrac{\omega}{T}) {\ut{$h$}}(\omega) \Big\}(x).
\end{align*}
Since $\Iop_{a,b}^{\Dord; \scalebox{0.4}{$k$}*}$ is scale-invariant of order $-\Dord$, we shall have
\begin{align*}
\varphi(x)-T^{k+1}\,\varphi(Tx) ~ \stackrel{\Iop_{a,b}^{\Dord; \scalebox{0.4}{$k$}*}}{\longmapsto} &~ \big(\Iop_{a,b}^{\Dord; \scalebox{0.4}{$k$}*}\,\varphi\big)(x) - T^{\scalebox{0.8}{\ut{$\Dord$}}+1} \big(\Iop_{a,b}^{\Dord; \scalebox{0.4}{$k$}*}\,\varphi\big)(Tx) \nonumber\\
&=  
\mathcal{F}^{-1}\Big\{  \big(\widehat{\phi}(\omega) - \widehat{\phi}(\tfrac{\omega}{T}) \big){\ut{$h$}}(\omega) \Big\}(x)\nonumber\\
&=  
\tfrac{1}{2\pi} \int_{\R} \big(\widehat{\phi}(\omega) - \widehat{\phi}(\tfrac{\omega}{T}) \big){\ut{$h$}}(\omega) \ee^{\jj\omega x}\dd \omega \nonumber\\
&=  
\tfrac{\jj}{2\pi \, x} \int_{\R} \tfrac{\dd}{\dd \omega} \Big(\big(\widehat{\phi}(\omega) - \widehat{\phi}(\tfrac{\omega}{T}) \big){\ut{$h$}}(\omega) \Big)\ee^{\jj\omega x}\dd \omega \nonumber\\
&=  \dots
= \tfrac{\jj^n}{2\pi \, x^n} \int_{\R} \tfrac{\dd^{n}}{\dd \omega^{n}} \Big(\big(\widehat{\phi}(\omega) - \widehat{\phi}(\tfrac{\omega}{T}) \big){\ut{$h$}}(\omega) \Big)\ee^{\jj\omega x}\dd \omega \nonumber\\
&=
\tfrac{\jj^n}{x^n} \mathcal{F}^{-1}\Big\{  \tfrac{\dd^n}{\dd \omega^n} \Big(\big(\widehat{\phi}(\omega) - \widehat{\phi}(\tfrac{\omega}{T}) \big){\ut{$h$}}(\omega) \Big) \Big\}(x),
\end{align*}
where $n=\lceil \real(\ut{$\Dord$})\rceil+1=k+1-\lfloor \real(\Dord)\rfloor$. 
Recalling Lemma \ref{lemma:L1property}, we have that $\|\tfrac{\dd^n}{\dd\omega^n} \big(\widehat{\phi}(\omega) - \widehat{\phi}(\tfrac{\omega}{T}) \big){\ut{$h$}}(\omega)\|_{1} 
\leq 
c\, T^{\real(\scalebox{0.8}{\ut{$\Dord$}}) +2-n} < c\, 4^{1-\{\real(\Dord)\}}$. Therefore,
\begin{align*}
\big|\big(\Iop_{a,b}^{\Dord; \scalebox{0.4}{$k$}*} \varphi\big)(x) - T^{\scalebox{0.8}{\ut{$\Dord$}}+1} \big(\Iop_{a,b}^{\Dord; \scalebox{0.4}{$k$}*} \varphi\big)(Tx)\big| 
 & \leq 
\tfrac{ \Big\|\mathcal{F}^{-1}\Big\{  \tfrac{\dd^n}{\dd \omega^n} \Big(\big(\widehat{\phi}(\omega) - \widehat{\phi}(\tfrac{\omega}{T}) \big){\ut{$h$}}(\omega) \Big) \Big\}(x) \Big\|_{\infty}}{ |x|^n }  \nonumber\\
& \leq 
\tfrac{ \Big\| \tfrac{\dd^n}{\dd \omega^n} \Big(\big(\widehat{\phi}(\omega) - \widehat{\phi}(\tfrac{\omega}{T}) \big){\ut{$h$}}(\omega) \Big)  \Big\|_{1}}{ |x|^n } \nonumber\\
& \leq
\tfrac{ c \, 4^{1-\{\real(\Dord)\} }}{|x|^n}.
\end{align*}
By setting $\Dord$ and $m$ parameters of Lemma \ref{lemma:ArashTechnique} as $\ut{$\Dord$}+1$ and $n$, respectively, we observe that
\begin{align}\label{eq:ProvdDecayRateIntegratorAdj}
\exists\,\bar{c}>0,~\forall\,1\leq|x|:~~~ \big|\big(\Iop_{a,b}^{\Dord; \scalebox{0.4}{$k$}*} \varphi\big)(x)\big| \leq \left\{\begin{array}{ll}
\tfrac{\bar{c}}{|x|^{k+1-\real(\Dord)}}, & \real(\Dord)\not\in \N, \phantom{\Big|}\\
\bar{c}\tfrac{1+\log|x|}{|x|^{k+1-\real(\Dord)}}, & \real(\Dord)\in \N. \phantom{\Big|}\\
\end{array}
\right.
\end{align}
Equation \eqref{eq:ProvdDecayRateIntegratorAdj} proves all the decay results stated in Theorem \ref{theo:Integrator1}. The above result further shows that
\begin{align} \label{eq:DecayProof}
\int_{|x|\geq 1}\big|\big(\Iop_{a,b}^{\Dord; \scalebox{0.4}{$k$}*} \varphi\big)(x)\big|^p\dd x <\infty,
\end{align}
for $p>\frac{1}{k+1-\real(\Dord)}$, whether $\real(\Dord)\in\N$ or $\real(\Dord)\not\in\N$.  Now, it is easy to combine \eqref{eq:SmallK0}, \eqref{eq:LargeK0} and \eqref{eq:DecayProof} to conclude that $\Iop_{a,b}^{\Dord; \scalebox{0.4}{$k$}*} \varphi$ is in $L_p$ for $p>\frac{1}{k+1-\real(\Dord)}$ if $k=\lfloor \real(\Dord)\rfloor$, and in $L_p$ for $\frac{1}{k+1-\real(\Dord)}<p<\frac{1}{k-\real(\Dord)}$ if $k>\lfloor \real(\Dord)\rfloor$.
$\hspace{\stretch{1}}\blacksquare$

 				
\subsection{Proof of Theorem \ref{theo:ImpulseResponse}}
As the operators are primarily defined in the Fourier domain, we use the fact that $\mathcal{F}\{\delta(\cdot-\tau)\}(\omega) = \ee^{-\jj\omega \tau}$. With this, the result regarding $\Dop_{a,b}^{\Dord}$ directly follows from \eqref{eq:InvF_Wp}:
\begin{align*}
\big(\Dop_{a,b}^{\Dord} \delta(\cdot-\tau)\big)(x) 
&= \mathcal{F}^{-1}\bigg\{ \ee^{-\jj\omega \tau} h_{a,b}^{\Dord}(\omega) \bigg\}(x) 
= \mathcal{F}^{-1}\bigg\{ h_{a,b}^{\Dord}(\omega) \bigg\}(x-\tau) \nonumber\\
&= a \mathcal{F}^{-1}\bigg\{ \omega_{+}^{\Dord} \bigg\}(x-\tau) + b \mathcal{F}^{-1}\bigg\{ \omega_{-}^{\Dord} \bigg\}(x-\tau) \nonumber\\
&=a \mathcal{F}^{-1}\bigg\{ \omega_{+}^{\Dord} \bigg\}(x-\tau) + b \mathcal{F}^{-1}\bigg\{ \omega_{+}^{\Dord} \bigg\}(\tau- x) \nonumber\\
&= \tfrac{\Gamma(\Dord+1)}{2\pi}\Big(\tfrac{a}{(\jj \tau - \jj x)^{\Dord+1}} + \tfrac{a}{(\jj x - \jj \tau)^{\Dord+1}}\Big).
\end{align*}

The claims for both $\Iop_{a,b}^{\Dord; \scalebox{0.6}{$k$}}$ and $\Iop_{a,b}^{(\Dord; \scalebox{0.6}{$k$})*}$ are obtained via Lemma \ref{lemma:InverseFourier}. We first show the result for $\Iop_{a,b}^{\Dord; \scalebox{0.6}{$k$}}$ and then, proceed with $\Iop_{a,b}^{(\Dord; \scalebox{0.6}{$k$})*}$:
\begin{align}\label{eq:ImpRespL}
&\big(\Iop_{a,b}^{\Dord; \scalebox{0.6}{$k$}} \delta(\cdot-\tau)\big)(x) 
= 
\tfrac{1}{2\pi} \int_{\R} \ee^{-\jj \omega \tau} \tfrac{\ee^{\jj \omega x} - \sum_{j=0}^{k-1}\tfrac{(\jj\omega x)^j}{j!}}{ h_{a,b}^{\Dord}(\omega) }\dd\omega \nonumber\\
&= 
\tfrac{1}{2a\pi} \int_{\R} \ee^{-\jj \omega \tau} w_{+}^{-\Dord} \Big(\ee^{\jj \omega x} \hspace{-1mm}-\hspace{-1mm} \sum_{j=0}^{k-1}\tfrac{(\jj\omega x)^j}{j!}\Big) \dd\omega 
+
\tfrac{1}{2b\pi} \int_{\R} \ee^{-\jj \omega \tau} w_{-}^{-\Dord} \Big(\ee^{\jj \omega x} \hspace{-1mm}-\hspace{-1mm} \sum_{j=0}^{k-1}\tfrac{(\jj\omega x)^j}{j!}\Big) \dd\omega\nonumber\\
&= 
\tfrac{1}{a}  \mathcal{F}^{-1}\bigg\{  w_{+}^{-\Dord} \Big(\ee^{\jj \omega x} \hspace{-1mm}-\hspace{-1mm} \sum_{j=0}^{k-1}\tfrac{(\jj\omega x)^j}{j!}\Big)   \bigg\}(-\tau) 
+
\tfrac{1}{b} \mathcal{F}^{-1}\bigg\{ w_{-}^{-\Dord} \Big(\ee^{\jj \omega x} \hspace{-1mm}-\hspace{-1mm} \sum_{j=0}^{k-1}\tfrac{(\jj\omega x)^j}{j!}\Big) \bigg\}(-\tau) \nonumber\\
&=
\tfrac{ \mathcal{F}^{-1}\Big\{  w_{+}^{\scalebox{0.6}{\ut{$\Dord$}}-k} \big(\ee^{\jj \omega x} - \sum_{j=0}^{k-1}\frac{(\jj\omega x)^j}{j!}\big)   \Big\}(-\tau) }{a}   
+
\tfrac{ \mathcal{F}^{-1}\Big\{ w_{+}^{\scalebox{0.6}{\ut{$\Dord$}}-k} \big(\ee^{-\jj \omega x} - \sum_{j=0}^{k-1}\tfrac{(-\jj\omega x)^j}{j!}\big) \Big\}(\tau) }{b} ,
\end{align}
where $\scalebox{1}{\ut{$\Dord$}}=k-\Dord$ and  $\real (\scalebox{1}{\ut{$\Dord$}}) >-1$. Now, the two inverse Fourier terms in \eqref{eq:ImpRespL} could be evaluated via Lemma \ref{lemma:InverseFourier}. The desired claim form then follows with simple modifications.

For operator $\Iop_{a,b}^{(\Dord; \scalebox{0.6}{$k$})*}$, we follow a similar approach:
\begin{align}\label{eq:ImpRespLAdj}
&\big(\Iop_{a,b}^{(\Dord; \scalebox{0.6}{$k$})*} \delta(\cdot-\tau)\big)(x) = \mathcal{F}^{-1}\bigg\{ \tfrac{\ee^{-\jj\omega \tau}-\sum_{j=0}^{k-1}\tfrac{(-\jj\omega \tau)^j}{j!} }{h_{a,b}^{\Dord}(-\omega)} \bigg\}(x) \nonumber\\
&= 
\tfrac{ \mathcal{F}^{-1}\Big\{ w_{+}^{-\Dord}\big(\ee^{-\jj\omega \tau}-\sum_{j=0}^{k-1}\tfrac{(-\jj\omega \tau)^j}{j!} \big) \Big\}(x) }{b} 
+ 
\tfrac{ \mathcal{F}^{-1}\Big\{ w_{-}^{-\Dord}\big(\ee^{-\jj\omega \tau}-\sum_{j=0}^{k-1}\tfrac{(-\jj\omega \tau)^j}{j!} \big) \Big\}(x) }{a} \nonumber\\
&= 
\tfrac{ \mathcal{F}^{-1}\Big\{ w_{+}^{\scalebox{0.6}{\ut{$\Dord$}}-k}\big(\ee^{-\jj\omega \tau}-\sum_{j=0}^{k-1}\tfrac{(-\jj\omega \tau)^j}{j!} \big) \Big\}(x) }{b} 
+ 
\tfrac{ \mathcal{F}^{-1}\Big\{ w_{+}^{\scalebox{0.6}{\ut{$\Dord$}}-k}\big(\ee^{\jj\omega \tau}-\sum_{j=0}^{k-1}\tfrac{(\jj\omega \tau)^j}{j!} \big) \Big\}(-x) }{a} .
\end{align}
Again, the two inverse Fourier terms in \eqref{eq:ImpRespLAdj} could be evaluated via Lemma \ref{lemma:InverseFourier}, and simplified to yield the claimed form.

\section{Conclusion}
\label{sec:conclusion}
In this paper, we initially studied complex-order fractional operators. We showed that by applying certain complex-order fractional integration operators to real-valued symmetric $\alpha$-stable white noise processes, we can generate self-similar stable processes with complex-valued Hurst exponent. Some of the introduced processes can be whitened by applying complex-order fractional derivative operators; as such, they can be described as solutions of fractional complex-order stochastic differential equations. 
  While we proved the existence of the random processes using characteristic functionals, we provided tools for numerically approximating such processes. We further studied the smoothness properties of the processes and showed that they have stationary increments of large enough orders.


\section*{Acknowledgement}
Julien Fageot is supported by the Swiss National Science Foundation (SNSF) under Grant P400P2\_194364. 


\bibliographystyle{elsarticle-num} 
\bibliography{refs.bib}

\end{document}